\documentclass[11pt]{article} 
\pdfoutput=1
\usepackage{amssymb} 
\usepackage{cite}
\usepackage{dsfont}
\usepackage{epsfig} 
\usepackage{rotate}
\usepackage{rotating}

\newtheorem{lemma}{Lemma}

\newcommand{\Q}{\mathbb{Q}}

\newcommand{\Z}{\mathbb{Z}}

\newcommand{\beq}{\begin{equation}} 
\newcommand{\eeq}{\end{equation}}
\newcommand{\Leq}[1]{\label{#1}\end{equation}}
\newcommand{\bdm}{\begin{displaymath}} 
\newcommand{\edm}{\end{displaymath}}

\title{Continuous self-similarity \\
in parametric piecewise isometries}
 \author{J. H. Lowenstein and F. Vivaldi\dag}
\date{
{\small\it
Dept.~of Physics, New York University, 2 Washington Place, New York, NY 10003, USA
\\
\dag School of Mathematical Sciences, Queen Mary, University of London, London E1 4NS, UK
}
}
\begin{document} 
\maketitle
\begin{abstract}
We exhibit two distinct renormalization scenarios in many-parameter families
of piecewise isometries (PWI) of a rhombus.
The rotational component, defined over the quadratic field 
$\mathbb{K}=\mathbb{Q}(\sqrt{5})$, is fixed.
The translations are specified by affine functions of the parameters, with coefficients 
in $\mathbb{K}$. In each case the parameters range over a convex domain.

In one scenario the PWI is self-similar if and only if one parameter belongs to $\mathbb{K}$, 
while the other is free.
Such a continuous self-similarity is due to the possibility of merging adjacent atoms of an induced PWI, 
a common phenomenon in the Rauzy-Veech induction for interval exchange transformations.

In the second scenario, the phase space splits into several disjoint (non-convex) invariant 
components. We show that each component has continuous self-similarity, but due to the transversality
of the corresponding foliations, full self-similarity in phase space is achieved if and only
if both parameters belong to $\mathbb{K}$. 

All our computations are exact, using algebraic numbers.
\end{abstract}

\centerline{\small\it\today}

\newpage

\section{Introduction}\label{section:Introduction}

This work represents a first study of renormalization in many-parameter families
of planar piecewise isometries (PWI).
These are maps of polygonal domains partitioned into convex sub-domains ---called
atoms--- in such a way that the restriction of the map
to each atom  is an isometry.  The first-return map to any convex sub-domain $D$ is
a new PWI, called the induced PWI on $D$. If by repeating the induction we obtain
a sub-system conjugate to the original one via a suitable group of isometries and 
homotheties, then we consider the original PWI to be renormalizable.

Recent work on renormalization in two-dimensional parametric families concerned 
one-parameter deformations of the translational part of a PWI.
\cite{Hooper,Schwartz:14,LowensteinVivaldi:14}. Induction is accompanied by a 
transformation $s\mapsto r(s)$, where $s$ is the parameter and $r$ is the renormalization
function. Self-similarity then corresponds to the periodic points of $r$.

This setting is analogous to Rauzy-Veech induction for interval exchange transformations
(IET's, see \cite{Rauzy,Veech,Yoccoz}), which are one-dimensional PWI's.
For an IET, the parameters are a vector of sub-interval lengths together with a permutation, 
and the renormalization acts on the lengths via an integral matrix. 
The fixed-point condition for self-similarity is an eigenvalue condition for a product of matrices.

An arithmetical characterisation of self-similarity is provided by the Boshernitzan-Carroll
theorem \cite{BoshernitzanCarroll}, which states that if an IET is defined over a quadratic
number field (meaning that all intervals' lengths belong to that field), then inducing on atoms
results in only finitely many distinct IETs, up to scaling.
In the case of two intervals (a rotation), this theorem reduces to Lagrange's theorem on 
the eventual periodicity of the continued fractions coefficients of quadratic irrationals. 
However, unlike for continued fractions,
there are self-similar IETs over fields of larger degree \cite{ArnouxYoccoz}.
It is also known that in a uniquely ergodic self-similar IET, the scaling constant is a unit
in a distinguished ring of algebraic integers \cite{PoggiaspallaLowensteinVivaldi}.

In two dimensions general results are scarce \cite{Poggiaspalla:03,Poggiaspalla:06}.
All early results on renormalization concerned specific models of PWI's
defined over quadratic fields (the field of a PWI is determined by the entries of the 
rotation matrices and the translation vectors defining the isometries)
\cite{LowensteinHatjispyrosVivaldi,AdlerKitchensTresser,KouptsovLowensteinVivaldi,%
AkiyamaBrunottePethoSteiner,Schwartz:09}. 
A more intricate form of renormalization has been found in a handful of cubic cases 
\cite{GoetzPoggiaspalla,LowensteinKouptsovVivaldi}. 

The first results on parametric families concerned polygon-exchange transformations,
due to Hooper \cite{Hooper} (on the measure of the periodic and aperiodic sets in 
a two-parameter family of rectangle-exchange transformations) and Schwartz 
\cite{Schwartz:14} (on the renormalization group of a one-parameter family of 
polygon-exchange transformations). Subsequently, the present authors \cite{LowensteinVivaldi:14}
studied two one-parameter families of piecewise isometries. Each family has a fixed rotational 
component defined over a quadratic field ($\mathbb{Q}(\sqrt{5})$ and 
$\mathbb{Q}(\sqrt{8})$, respectively), and parameter-dependent translations.
It was shown that self-similarity occurs if and only if the parameter belongs to
the respective field.

\begin{figure}[h]
\hfil\epsfig{file=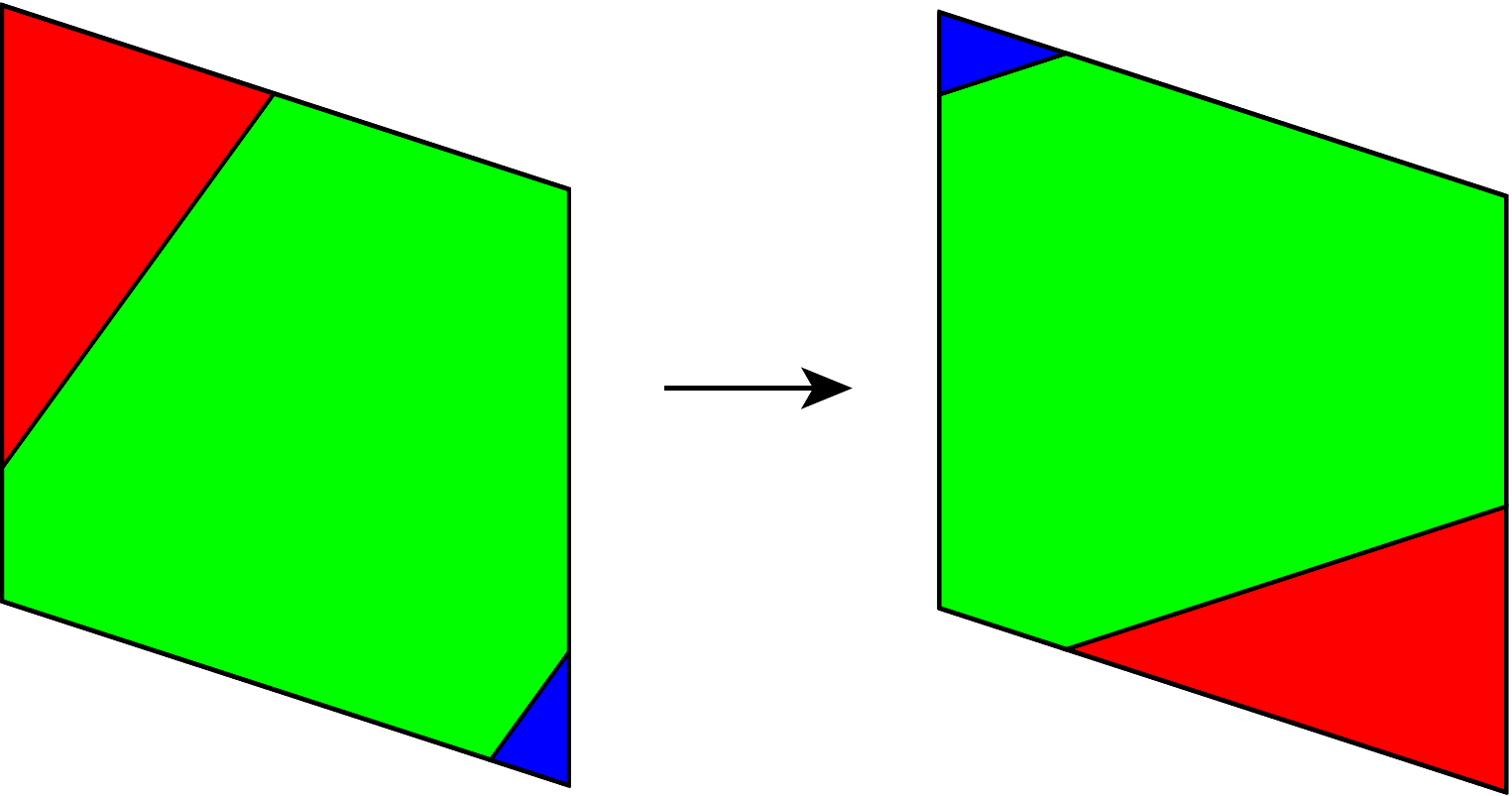,width=8cm}\hfil

\caption{\label{fig:3Rpwi}\small One-parameter rhombus map of the pentagonal model 
of \cite{LowensteinVivaldi:14}. The rhombus and all directions are fixed.
The parameter shifts the inner boundaries of the atoms in such a way that all
three atoms retain a reflection symmetry.
}
\end{figure}

The mapping for the pentagonal model is shown in figure \ref{fig:3Rpwi}.
For the parameter $s$ restricted to a suitable interval,
there is an induced PWI on a triangular sub-domain (the so-called {\it base triangle}), 
which reproduces itself after scaling and the reparametrisation $s\mapsto r(s)$. 
After an affine change of parameter, the function $r$ was found to be of L\"{u}roth type 
---a piecewise affine version of Gauss's map \cite{BarrionuevoEtAl,Galambos,Luroth}.  
In the pentagonal model, the discontinuities of $r$ accumulate at the origin; 
in the octagonal case one has $r=f\circ f$, where $f$ has two accumulation 
points of discontinuities. 
In both cases $r$ is expanding and preserves the Lebesgue measure.

\begin{figure}[h]
\hfil\epsfig{file=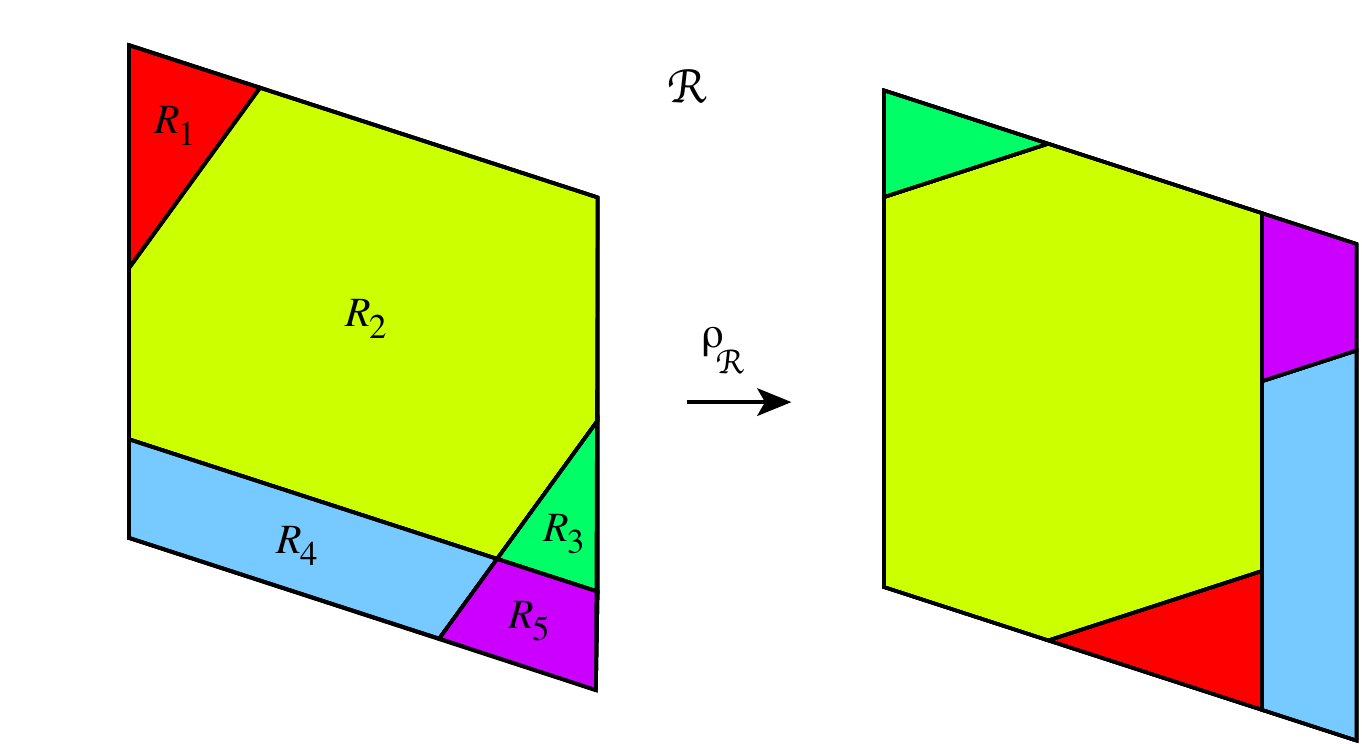,width=12cm}\hfil

\caption{\label{fig:5Rpwi}\small Two-parameter rhombus map.
}
\end{figure}
\begin{figure}[h]
\hfil\epsfig{file=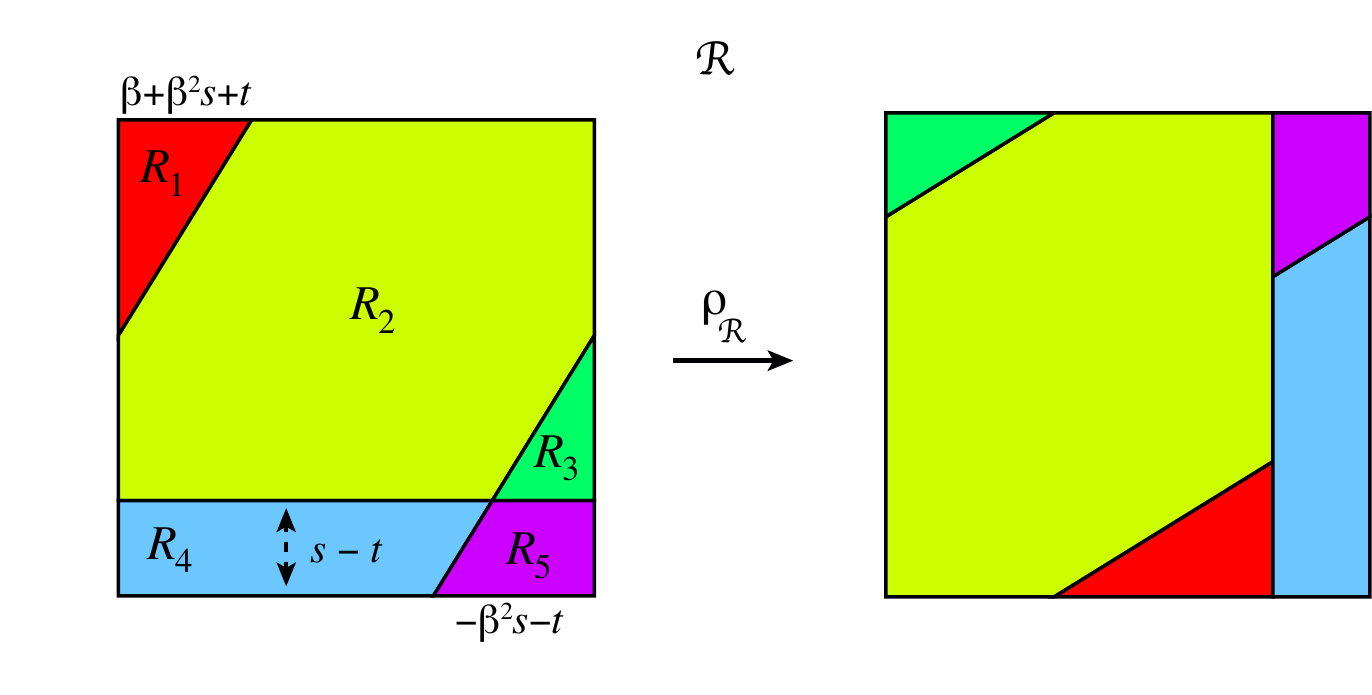,width=12cm}\hfil

\caption{\label{fig:5square}\small The two-parameter map conjugate to the
rhombus map of figure \ref{fig:5Rpwi}.
}
\end{figure}

The present work deals with two-parameter extensions of the one-parameter 
pentagonal model of \cite{LowensteinVivaldi:14} (figure \ref{fig:5Rpwi}).
The roles of the parameters is made explicit in figure \ref{fig:5square}, 
where we employ a new co-ordinate system in order to restrict the arithmetic
to the quadratic field $\mathbb{K}=\mathbb{Q}(\sqrt{5})$, rather than a 
bi-quadratic field ---see section \ref{section:Computations}.
In the new co-ordinates, the atom boundary lines are of the form 
$$
\mathbf{u}_i\cdot \mathbf{x}=b_i,
\qquad
b_i=b_{i,0}+ b_{i,1} s_1 + b_{i,2} s_2,\quad b_{i,j}\in \mathbb{Q}(\sqrt{5})
\qquad 0\leqslant i\leqslant 4
$$ 
where $s_1,s_2$ are the parameters.
We use these co-ordinates for calculations, reserving the original co-ordinates for graphics.

In section \ref{section:Preliminaries} below, we develop the necessary geometrical constructs:
tiles, dressed domains, isometries, etc.
In section \ref{section:BaseTriangle} we review the renormalization theory of
the base triangle developed in \cite{LowensteinVivaldi:14}, with two additions.
We extend the isometries to the atoms' boundaries (which will be needed to glue 
together atoms in section \ref{section:MainResults}), and we introduce an improved 
renormalization scheme whose renormalization function has finitely many singularities.

The first renormalization scenario is established in section \ref{section:MainResults}.  
There we make three successive inductions on triangular sub-domains, to obtain a PWI 
whose eight distinct return orbits tile the rhombus, apart from the orbits of finitely
many periodic domains.  
We identify a convex polygon $\Pi$ in the $t,s$ parameter space
within which the induced map has only three distinct isometries. 
By merging neighbouring atoms which share the same isometry, we recover the 
one-parameter PWI discussed in the previous section.
As a result, self-similar dynamics occurs if and only if 
$s\in\mathbb{K}$, and since $t$ is arbitrary, the self-similarity
constraint corresponds to a foliation of $\Pi$.

The appearance of a hidden reduction of the number of atoms, which reveals itself 
only after induction, is present in a simpler form in the Rauzy-Veech 
induction of interval exchange transformations. 
In section \ref{section:Rauzy} we show that the occurrence of free parameters
in renormalizable IETs can be understood in terms of the properties of the associated
translation surfaces \cite{KontsevichZorich}. In particular, there are self-similar
IETs without free parameters, as long as the number of intervals is greater than three.

In view of this analogy, we conducted an extensive search for a non-degenerate 
two-parameter family, using two-dimensional sections of a three-parameter system 
of the $2\pi/5$ rhombus (section \ref{section:Explorations}).
We found a weak form of non-degeneracy, resulting from the co-existence of two 
continuously renormalizable families with transversally intersecting foliations.
Each family corresponds to a (non-convex) invariant set in phase space, and hence to an 
ergodic component of the exceptional set.

In many-parameter families of piecewise isometries over quadratic fields, the existence of
an irreducible component in phase space which is self-similar without free parameters has 
yet to be established.
\bigskip

\noindent {\sc Acknowledgements:} \/
JHL and FV would like to thank, respectively, the School of
Mathematical Sciences at Queen Mary, University of London,
and the Department of Physics of New York University, for their hospitality.

\section{Preliminaries}\label{section:Preliminaries}
Throughout this paper, we let 
\begin{equation}\label{eq:AlphaOmega}
\alpha=\sqrt{5},\qquad \omega=(\alpha+1)/2,\qquad\beta=\omega^{-1}=(\alpha-1)/2.
\end{equation}
The arithmetical environment is the quadratic field $\mathbb{Q}(\omega)$ 
with its ring of integers $\mathbb{Z}[\omega]$, given by
\begin{equation}\label{eq:FieldRing}
\mathbb{Q}(\omega)=\{x+y\alpha\,:\,x,y\in\mathbb{Q}\},
\qquad
\mathbb{Z}[\omega]=\{m+n\omega\,:\,m,n\in\mathbb{Z}\}.
\end{equation}
The number $\omega$, which is the fundamental unit in $\Z[\omega]$
(see \cite[chapter 6]{Cohn}), will determine the scaling under
renormalization. The number $\beta=\omega-1$ is also a unit.

\subsection{Planar objects\label{section:PlanarObjects}}

A \textit{tile} $X$ with $n$ edges is a convex polygon defined by the half-plane
conditions
\begin{equation}
\begin{array}{lcll} 
\mathbf{u}_{m_i}\cdot\mathbf{x} &<& b_i & \mbox{(excluded edge)}\nonumber\\
&\mbox{or}&& \label{eq:HalfPlaneConditions}\\
\mathbf{u}_{m_i}\cdot\mathbf{x} &\geqslant& b_i & \mbox{(included edge)}\nonumber
\end{array}\qquad i=1,\ldots,n, 
\end{equation}
where $\mathbf{x}=(x,y)$, $b_i\in\mathbb{R}$, and the $\mathbf{u}_m$ are the
vectors
\beq\label{eq:uvector}
\mathbf{u}_m=\left(\cos\frac{2\pi m}{5},\sin\frac{2\pi m}{5}\right)\qquad m\in\{0,\ldots,4\}.
\eeq
For the $i$th edge, defined by $\mathbf{u}_{m_i}\cdot\mathbf{x}=b_i$,
we introduce an index $\epsilon_i$, where $\epsilon_i=-1$ if the edge is 
included in $X$, and $\epsilon_i=1$ if it is excluded.
We then represent $X$ as a triple of $n$-vectors
\begin{equation}\label{eq:Tile}
X=[(m_1,\ldots,m_n), (\epsilon_1,\ldots,\epsilon_n), (b_1,\ldots,b_n)].
\end{equation}
We shall assume that $n$ is minimal, namely that $X$ is not definable by
fewer conditions.

A \textit{tiling} ${\bf X}$ is a set of disjoint tiles,
$$
{\bf X}= \{X_1,\ldots,X_N\}
$$
and is associated with a domain ${\rm X}$ (union of tiles)
$$
{\rm X}=\bigcup_{k=1}^N X_k.
$$
Note that a domain need not be convex, or even connected.
Note further that thanks to (\ref{eq:HalfPlaneConditions}), if a pair of
tiles have disjoint interiors but share a common boundary segment, that
segment belongs to one and only one tile of the pair.
This allows the possibility of gluing together adjacent tiles without
disturbing the inclusion relation of the respective edges.

\subsection{Similarity group}
The transformation properties of planar objects are provided by a group 
$\mathfrak{G}$ which comprises the rotations and reflections of the symmetry 
group of the regular pentagon (the dihedral group $D_5$) together with 
translations in $\mathbb{K}^2$ and real scale transformations. 

We adopt the following notation:
\begin{itemize}
\item [] $\mathtt{U}_m$: reflection about the line generated by $\mathbf{u}_m$.
\item [] $\mathtt{R}_m$: rotation by the angle $2 m\pi/5$.
\item [] $\mathtt{T}_{\bf d}$: translation by ${\bf d}\in \mathbb{K}^2$.
\item [] $\mathtt{S}_{\eta}$: scaling by $\eta\in\mathbb{R}_+$.
\end{itemize}
We write $\mathcal{X}\sim\mathcal{Y}$ to indicate that $\mathcal{X}$ is
similar to $\mathcal{Y}$, i.e., that
$\mathcal{X}=\mathtt{G}(\mathcal{Y})$ for some $\mathtt{G}\in\mathfrak{G}$.  
As $\mathfrak{G}$ is a group, this is an equivalence relation.  
Within $\mathfrak{G}$ we distinguish two important subgroups: the
\textit{isometry group} $\mathfrak{I}$ generated by rotations, reflections, 
and translations, and the \textit{dynamical group} $\mathfrak{I}_+$,
generated by rotations and translations.

\subsection{Dressed domains and sub-domains}
A \textit{dressed domain} is a triple
\begin{equation}\label{eq:DressedDomain}
\mathcal{X}=({\mathrm{X}, \mathbf{X}, \rho}),
\end{equation}
where $\mathbf{X}=\{X_1,\ldots,X_N\}$ is a tiling of the domain $X$,
and $\rho=\{\rho_1,\ldots,\rho_N\}$, where $\rho_i\in\mathfrak{J}_+$ is an 
orientation-preserving isometry 
acting on the tile $X_k$.
Under the action of $\mathtt{G} \in \mathfrak{G}$, a dressed domain $\mathcal{X}$ 
transforms as
$$
\mathtt{G}(\mathcal{X})= \mathtt{G}({\mathrm{X}, \mathbf{X}, \rho}) =
(\mathtt{G}(\mathrm{X}),\{\mathtt{G}(X_1),\ldots,\mathtt{G}(X_k)\}, \mathtt{G}\circ\rho\circ \mathtt{G}^{-1})
$$
where the conjugacy acts componentwise.
To emphasize the association of a mapping $\rho$ with a particular dressed domain $\mathcal X$, 
we  use the notation $\rho_{\stackrel{}{\mathcal{X}}}$. 

Let $\mathcal{X}= (\mathrm{X}, \mathbf{X}, \rho_{\stackrel{}{\mathcal{X}}})$ be a dressed domain, 
and let $\mathrm{Y}$ be a sub-domain of $\mathrm{X}$.  
We denote by $\rho_{\stackrel{}{\mathcal{Y}}}$ the first-return map on $\mathrm{Y}$ induced by
$\rho_{\stackrel{}{\mathcal{X}}}$.  
We call the resulting dressed domain 
$\mathcal{Y}= (\mathrm{Y}, \mathbf{Y}, \rho_{\stackrel{}{\mathcal{Y}}})$ 
a \textit{dressed sub-domain} of $\mathcal{X}$, and write
\begin{equation}\label{eq:YsubX}
\mathcal{Y}\triangleleft\mathcal{X}.
\end{equation}

The dressed sub-domain relation (\ref{eq:YsubX}) is \textit{scale invariant}, namely 
invariant under an homothety. Indeed, if
$\mathtt{S}_\eta$ denotes scaling by a factor $\eta$, then in the data 
(\ref{eq:Tile}) specifying a tile, the orientations $m_k$ remain unchanged,
while the pentagonal coordinates $b_k$ scale by $\eta$. 
Moreover, the identity
$$
\mathtt{S}_\eta \mathtt{T}_\mathbf{d}\mathtt{R}_n= \mathtt{T}_{\eta \mathbf{d}}\mathtt{R}_n\mathtt{S}_\eta.
$$
shows that the piecewise isometries $\rho$ scale in the same way.
We conclude that the subdomain relation (\ref{eq:YsubX}) is preserved if the dressed 
domain parameters are scaled by the same factor for both members.

\subsection{Parametric dressed domains}\label{section:ParametricDressedDomains}
In this article we consider continuously deformable dressed domains
$\mathcal{X}=\mathcal{X}(s)$
called \textit{parametric dressed domains}, depending on
a real parameter vector $s=(s_1,\ldots,s_p)$.

These are domains whose tiles $X_k$ and image tiles $\rho_k(X_k)$ 
depend on $s$ only via the coefficients $b_i$, while the parameters 
$n$, $m_i$ and $\epsilon_i$ remain fixed [see (\ref{eq:Tile})].
We shall require that the $b_i$'s be affine functions of $s_1,\dots,s_p$, with coefficients
in $\mathbb{Q}(\omega)$. 
Algebraically, this is expressed as
\begin{equation}\label{eq:ModuleS}
b_i\in\mathbb{S}
\hskip 40pt 
{\mathbb{S}}=\Q(\omega) + \Q(\omega) s_1 +\cdots + \Q(\omega) s_p,
\end{equation}
where $s_1,\ldots,s_p$ are regarded as indeterminates. The set $\mathbb{S}$ is
is a ($p+1$)-dimensional vector space over $\Q(\omega)$ (a $\Q(\omega)$-module).

The condition (\ref{eq:ModuleS}) gives us affine functions 
$b_i:\mathbb{R}^p\to\mathbb{R}$
\begin{equation}\label{eq:b_i}
b_i(s_1,\ldots,s_p)=b_{i,0}+b_{i,1}s_1+\cdots+b_{i,p}s_p\qquad b_{i,j}\in \mathbb{Q}({\omega}).
\end{equation}
We define the \textit{bifurcation-free set} $\mathbf{\Pi}(\mathcal{X})$ to be the maximal 
open set such that all of the edges of all $X_k(s)$ have non-zero lengths. 
Note that other types of bifurcations may occur if $\mathcal{X}$ is 
embedded within a larger domain (see section \ref{section:Computations}.)


\subsection{Renormalizable dressed domains}\label{section:Renormalizability}

A parametric dressed domain $\mathcal{X}(s)$ is said to be \textit{renormalizable} 
over an open domain $\Pi\subset \mathbb{R}^p$ if there exists a piecewise smooth map 
$r:\Pi\rightarrow \overline \Pi$ such that for every choice of $s\in r^{-1}(\Pi)$ the 
dressed domain $\mathcal{X}(s)$ has a dressed subdomain $\mathcal{Y}$ similar to 
$\mathcal{X}(r(s))$ which satisfies the recursive tiling property.
The function $r$ depends only on $s$, a requirement of scale invariance. 
In general, we have $\mathcal{Y}=\mathcal{Y}_{i(s)}(s)$, where $i$ is a discrete index. 
The set $r^{-1}(\Pi)$ need not be connected (even if $i$ is constant), each connected
component being a bifurcation-free domain of $\mathcal{Y}$. 
(To extend the renormalization function $r$ to the closure of $\Pi$, one must include
bifurcation parameter values, as in \cite{LowensteinVivaldi:14}.)

If $s=s_0$ is eventually periodic under $r$, then we say that $\mathcal{X}(s_0)$ 
is \textit{self-similar}. A self-similar system has an induced sub-system which reproduces itself
on a smaller scale under induction. 

Let a parametric dressed domain $\mathcal{X}(s)$ have induced $\mathcal{X}_j(s)$, 
such that, for $j=1,\ldots,N$ we have: ($i$) $\mathcal{X}_j$ is renormalizable over a domain $\Pi_j$;
($ii$) the $\mathcal{X}_j$ recursively tile $\mathcal{X}$; ($iii$) the
$\Pi_j$ have non-empty intersection $\Pi$. Then we still consider $\mathcal{X}$ 
renormalizable over $\Pi$.

The definition of renormalizability given above is tailored to our model; 
it is not the most general possible, and it is local in parameter space. 
We allow $\mathcal{Y}$ to depend on a discrete
index (as in Rauzy induction for interval-exchange transformations
---see section \ref{section:Rauzy}) to obtain a simpler renormalization
function $r$ (section \ref{section:BaseTriangle}). 
We only require $\mathcal{X}$ to be eventually renormalizable, and we
allow $\mathcal{X}$ to have sub-domains with independent renormalization schemes 
(which is a common phenomenon, see section \ref{section:Explorations}).

\subsection{Computations}\label{section:Computations}
For computations, we use the 
Mathematica\textsuperscript\textregistered
procedures listed in
the Electronic Supplement \cite{ESupplement}.
All computations reported in this work are exact, employing integer 
and polynomial arithmetic, and the symbolic representation of algebraic
numbers.

The geometrical objects defined in section \ref{section:PlanarObjects} require 
arithmetic in a bi-quadratic field, since only the first component of 
the vectors $\mathbf{u}_m$ is in $\mathbb{Q}(\omega)$.
To circumvent this difficulty, we conjugate our PWI to a map of a square where the 
clockwise rotation $2\pi/5$ is represented by the following matrix over $\Z[\omega]$
$$
\left(\begin{array}{cc}0&1\\-1&\beta\end{array}\right)
$$
where $\beta$ was defined in (\ref{eq:AlphaOmega}).
(This is still a PWI with respect to a non-Euclidean metric.)
In the new co-ordinates, the vectors $\mathbf{u}_m$ become 
$$
\{(1,0),(0,1),(-1,\beta),(-\beta,-\beta),(\beta,-1)\}
$$
which belong to $\mathbb{Z}[\omega]^2$.
With this representation, all of our calculations can be performed within 
the module $\mathbb{S}$ defined in (\ref{eq:ModuleS}).
We shall still display our figures in the original coordinates, where 
geometric relations (especially reflection symmetries) are more apparent 
to the eye.

In constructing a return map orbit of a domain $\mathcal{X}(s)$ by direct iteration, one 
determines inclusion and disjointness relations among domains, which require evaluations 
of inequalities (\ref{eq:HalfPlaneConditions}).
Since the latter are expressed by affine functions of the parameter $s$ in some polytope 
$\Pi$, it suffices to check the inequalities on the boundary of $\Pi$. 
All these boundary values belong to the field $\Q(\omega)$, and the inequalities can be
reduced to integer inequalities. 

Typically, $\mathcal{X}$ will be immersed in a larger domain $\mathcal{Y}$ (an atom, say). 
Therefore, in addition to the intrinsic bifurcation-free polytope $\Pi(\mathcal{X})$ defined 
in section \ref{section:ParametricDressedDomains}, one must also consider the polytope 
$\Pi(\mathcal{X},\mathcal{Y})$ defined by the inclusion $\mathcal{X}(s)\subset\mathcal{Y}(s)$, 
as well as intersection of these polytopes.

The recursive tiling property defined in section \ref{section:Renormalizability} is established
by adding up the areas of the tiles of all the orbits, and comparing it with the total area 
of the parent domain. 

With these techniques, we are able to establish rigorously a variety of statements valid 
over convex sets in parameter space.

\section{Base triangle}\label{section:BaseTriangle}
The \textit{base triangle} is the simplest one-parameter renormalizable 
piecewise isometry associated with rotations by $2\pi/5$; it is 
self-similar precisely for parameter in the quadratic field $\mathbb{Q}(\sqrt{5})$.
It was instrumental to the proof of renormalizability of a one-parameter rhombus 
map in \cite{LowensteinVivaldi:14}, and it will appear again in the 
many-parameter versions presented here.

We develop a variant of the model presented in \cite{LowensteinVivaldi:14}, 
which includes boundary segments in the tiles and an improved renormalization scheme. 
The base triangle $\mathcal{B}$ prototype is the following dressed domain 
(see figure \ref{fig:Bpwi}):
\begin{figure}[h]
\hfil\epsfig{file=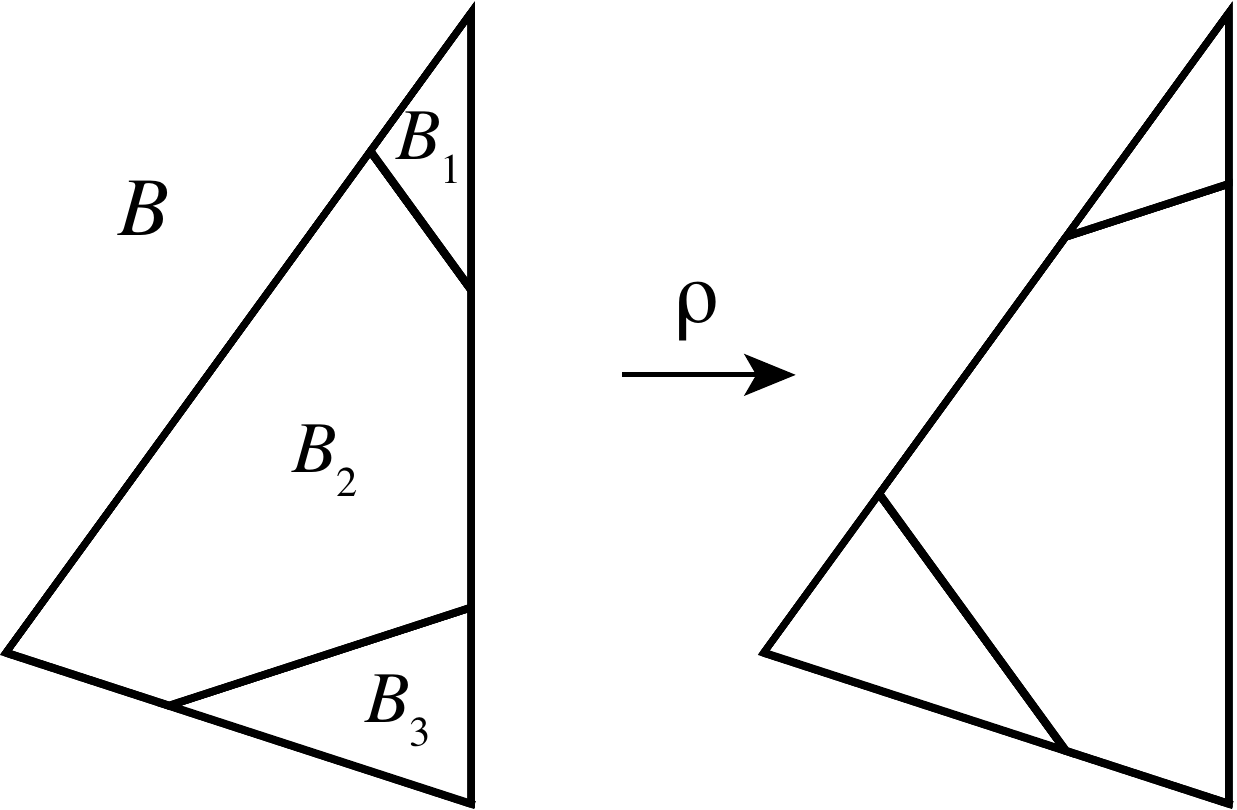,width=8 cm}\hfil
\caption{\label{fig:Bpwi} \small Base triangle prototype.} 
\end{figure}
$$
\mathcal{B}=(B, (B_1,B_2,B_3), (\rho_1,\rho_2,\rho_3))
$$
where
\begin{eqnarray}
B&=&[(1,0,2),(-1,1,1),(\tau-\omega^2,0,0)],\nonumber\\
B_1&=&[(0,2,3),(1,1,1),(0,0,\omega-\omega\tau)],\label{eq:Bdef}\\
B_2&=&[(1,4,0,3,2),(-1,1,1,-1,1),(\tau-\omega^2,\omega^2-\omega \tau,0, \omega-\omega \tau,0)],\nonumber\\\
B_3&=&[(1,0,4),(-1,1,-1),(\tau-\omega^2,0,\omega^2-\omega \tau)].\nonumber
\end{eqnarray}
The dynamics is given by a local reflection of each atom about its own
symmetry axis, followed by a global reflection about the symmetry 
axis of $B$, which can be written as:
\begin{eqnarray}
\rho_1&=&\mathtt{T}_{(\omega\tau-\omega,-\omega^2+\omega^2 \tau)}\, \mathtt{R}_2\nonumber\\
\rho_2&=&\mathtt{T}_{(0,\omega^2 \tau-\omega^3)}\, \mathtt{R}_3 \label{eq:rhoBdef}\\
\rho_3&=&\mathtt{T}_{(\omega \tau-2\omega,\omega^2\tau-2\omega^2)}\, \mathtt{R}_2.\nonumber
\end{eqnarray}
Here we have chosen a coordinate system such that the peak of the isosceles triangle is at the origin 
and the altitude of the atom $B_3$ is the parameter $\tau$, which varies over the interval $(0,1)$ 
without the occurrence of a bifurcation.   
This parameter (together with time-reversal invariance) determines the scale-invariant properties 
of the dressed domain, since it is related to the ratio $\eta$ of altitudes of $B_3$ and $B$ by the formula
$$
\eta=\frac{\tau}{\omega^2-\tau}.
$$
As $\tau$ varies from $0$ to $1$, $\eta$ increases from $0$ to $\beta$. 

The edges of the domain $\mathcal{B}$ are included or excluded as stipulated in 
section \ref{section:Preliminaries}; a vertex joining two included edges is included, but is excluded otherwise.  
The renormalizability analysis will also require a second base triangle $\tilde{\mathcal{B}}$, 
differing from $\mathcal{B}$ by a change of sign of all edge coordinates and translation vectors, 
as well as of the respective $\epsilon_i$.   
The dressed domains $\mathcal{B}$ and $\tilde{\mathcal{B}}$ are $\mathfrak{G}$-inequivalent: not 
only do they have different boundary conditions, but their interiors differ by a rotation by 
$\pi$, not an element of the similarity group.

The renormalizability analysis for the base triangle is summarized in the following lemma:

\begin{lemma} \label{lemma:BaseTriangle}
Let $\mathcal{B}$ be as above. The following holds:
\begin{enumerate}
\item[\textrm{(i)}] For $0<\tau< \beta^2$, $\mathcal{B}$ has a dressed subdomain 
$\mathcal{B}_1\sim \mathcal{B}$ which is scaled by a factor $(1-\tau)/(\omega^2-\tau)$ and has shape 
parameter $r(\tau)=\omega^2 \tau$. 

\item[\textrm{(ii)}] For $\beta^2 < \tau <\beta$, $\mathcal{B}$ has a dressed subdomain 
$\mathcal{B}_2\sim \tilde{\mathcal{B}} $ which is scaled by a factor $\tau/(\omega^2-\tau)$ and 
has shape parameter $r(\tau)=\omega^3 (\beta-\tau)$. 

\item[\textrm{(iii)}] For $\beta < \tau <1$, $\mathcal{B}$ has a dressed subdomain 
$\mathcal{B}_3\sim \tilde{\mathcal{B}}$ which is scaled by a factor $\tau/(\omega^2-\tau)$ 
and has shape parameter $r(\tau)=\omega^2 (1-\tau)$. 
\end{enumerate}
\end{lemma}
 
\begin{figure}[h]
\hfil\epsfig{file=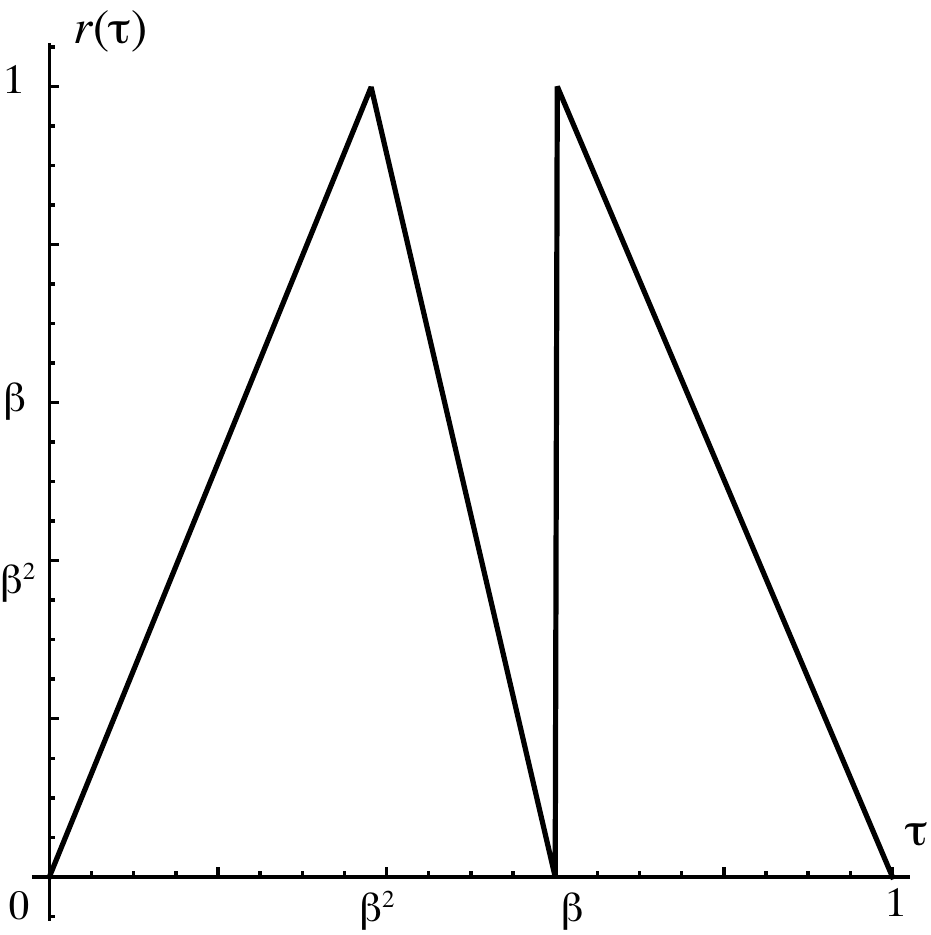,width=7 cm}\hfil
\caption{\label{fig:rplot} \small Renormalization function $r(\tau)$ for base triangles.} 
\end{figure}

The renormalization function $r$ has three branches (see figure \ref{fig:rplot}). 
In cases (ii) and (iii) one induces on the atom $B_3$, over two disjoint bifurcation-free parameter
ranges. 
Since the size of $B_3$ vanishes as $\tau$ approaches $0$, in the range (i) we induce on the triangle 
$[(1,0,2),(-1,1,1),(\beta^2 \tau-1,0,0)]$, which is not an atom.
This device prevents the occurrence of infinitely many singularities in the renormalization 
function found in \cite{LowensteinVivaldi:14}.

\begin{figure}[h!]
\hfil\epsfig{file=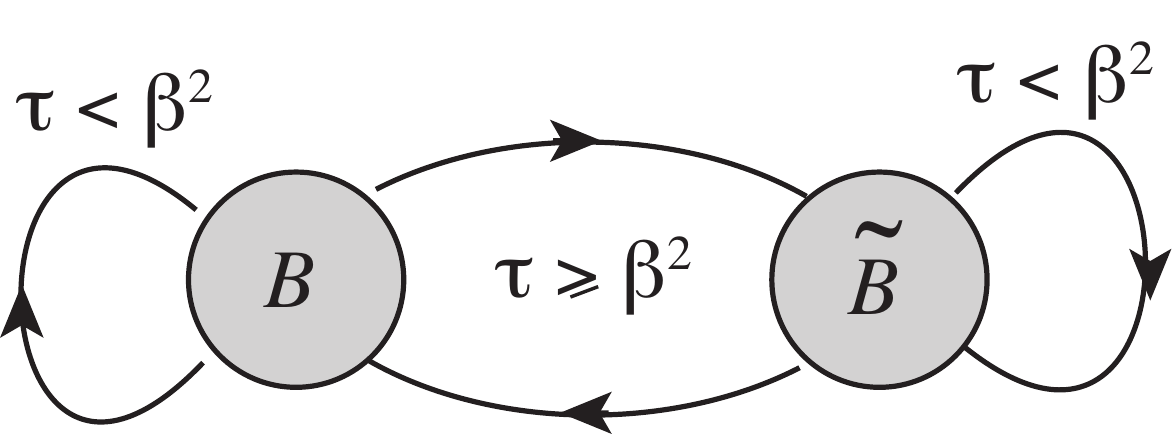,width=7 cm}\hfil
\caption{\label{fig:Bgraph} \small Renormalization graph for base triangles. 
A directed link from $X$ to $Y$ indicates that $X$ has an induced dressed subdomain equivalent to $Y$, 
subject to the the indicated parameter constraint.} 
\end{figure}

In each case, the return orbits of the atoms of $\mathcal{B}_i$, together with a finite number 
of periodic tiles, completely tile the triangle $B$.
For $\tilde{\mathcal{B}}$, the prescriptions (i)-(iii) hold with the roles of $\mathcal{B}$ and 
$\tilde{\mathcal{B}}$ exchanged. 
The induction relations are represented as the graph in figure \ref{fig:Bgraph}. 

As in \cite{LowensteinVivaldi:14}, the proof of Lemma 1 is by direct iteration, as discussed 
in section \ref{section:Computations}.  
The main difference in the two computational algorithms lies in the procedures used to verify 
inclusion and disjointness of tiles (see \cite{ESupplement}). 
Specifically, checking the sub-polygon relation $X\subset Y$ requires  verifying that no included 
vertex of $X$ has landed on an excluded edge of $Y$. 
Similarly, to decide that $X$ and $Y$ are disjoint, one must check that no included vertex of either 
tile lies on an included edge of the other.

With reference to lemma \ref{lemma:BaseTriangle}, we remark that the base triange $\mathcal{B}$, 
with its definition extended to the two-atom limiting cases $\tau=0, 1$, is in fact renormalizable 
also at the parameter values $0, \beta^2, \beta, 1$, with $r(\tau)=0$ in all these cases 
(see \cite{ESupplement} for the calculations).  
These additional parameter values are needed to make $\mathcal{B}$ renormalizable over 
the whole interval $[0,1]\cap \mathbb{Q}(\omega)$.  
We shall use this property in sections \ref{section:MainResults} and \ref{section:Explorations}.

\section{Continuous self-similarity}\label{section:MainResults}
We now turn to the two-parameter rhombus map introduced in section \ref{section:Introduction}
---see figures \ref{fig:5Rpwi} and \ref{fig:5square}. 
In suitable coordinates, the dressed domain is given by
$$
\mathcal{R}=(R,(R_1,\ldots,R_5), (\rho_{R_1},\ldots,\rho_{R_5}) ),
$$
with (see figure \ref{fig:5Rpwi})
\begin{eqnarray}
R&=& [(0, 1, 0, 1), (-1, -1, 1, 1), (-t, -s, 1 - t, 1 - s)],\nonumber\\
R_1&=&[(0, 2, 1), (-1, -1, 1), (-t, -s, 1 - s)],\nonumber\\
R_2&=&[ (0, 1, 2, 0, 1, 2), (-1, -1, -1, 1, 1, 1), (-t, -t, -1 - s,
   1 - t, 1 - s, -s)],\nonumber\\
R_3&=&[ (0, 2, 1), (1, 1, -1), (1 - t, -1 - s, -t)],\label{eq:Rdef}\\
R_4&=&[ (0, 1, 2, 1), (-1, -1, -1, 1), (-t, -s, -1 - s, -t)],\nonumber\\
R_5&=&[ (1, 0, 1, 2), (-1, 1, 1, 1), (-s, 1 - t, -t, -1 - s)]\nonumber
\end{eqnarray}
\beq\label{eq:rhoRdef}
\rho_{R_1}=\mathtt{T}_{(0,0)}\, \mathtt{R}_4,\qquad
\rho_{R_2}=\mathtt{T}_{(0,1)}\, \mathtt{R}_4,\qquad
\rho_{R_3}=\mathtt{T}_{(0,2)}\, \mathtt{R}_4,\qquad
\eeq
$$
\rho_{R_4}=\mathtt{T}_{(1,1)}\, \mathtt{R}_4.\qquad
\rho_{R_5}=\mathtt{T}_{(1,2)}\, \mathtt{R}_4.
$$

The corresponding bifurcation-free parametric domain $\Pi(\mathcal{R})$,
defined in section \ref{section:ParametricDressedDomains},
is found to be the triangle with vertices at 
$(0,0)$, $(-1/\alpha,-1/\alpha)$, and $(\beta/\alpha,-\beta/\alpha)$.
On the boundary of $\Pi(\mathcal{R})$ given by with $s=t$, the dressed domain $\mathcal{R}$
collapses into the one-parameter pentagonal model of \cite{LowensteinVivaldi:14}, 
and hence is self-similar for all $s\in\Q(\omega)$ within a suitable interval.

\begin{figure}[h]
\hfil\epsfig{file=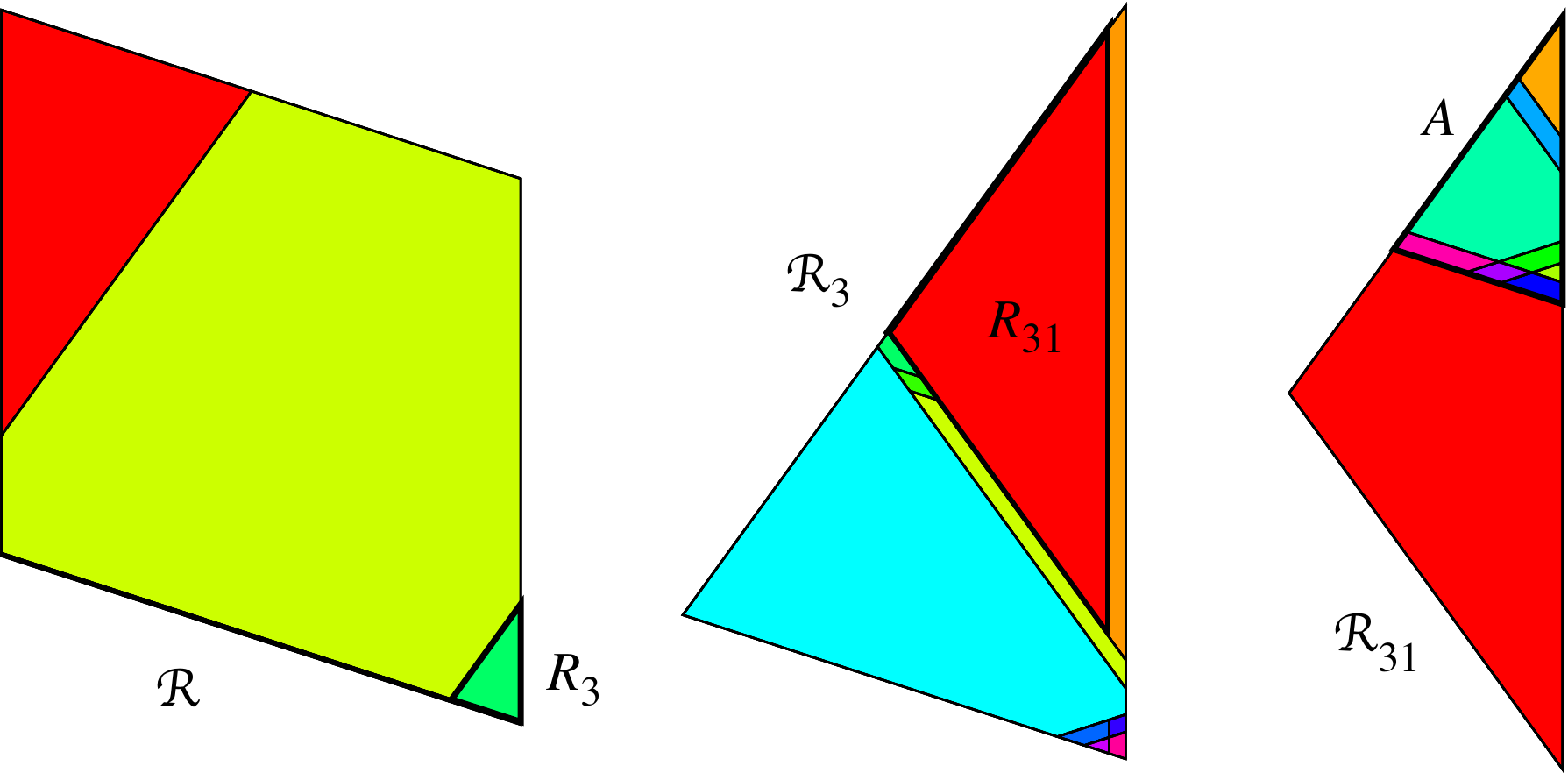,width=12cm}\hfil
\caption{\label{fig:3Xinduction}\small 
The first two steps of the triple induction $\mathcal{R}\triangleright\mathcal{R}_3\triangleright\mathcal{R}_{31}\triangleright\mathcal{A}$.  
The third step produces the dressed domain $\mathcal{A}$ shown in figure \ref{fig:Apwi}.
}
\end{figure}

\begin{figure}[h]
\hfil\epsfig{file=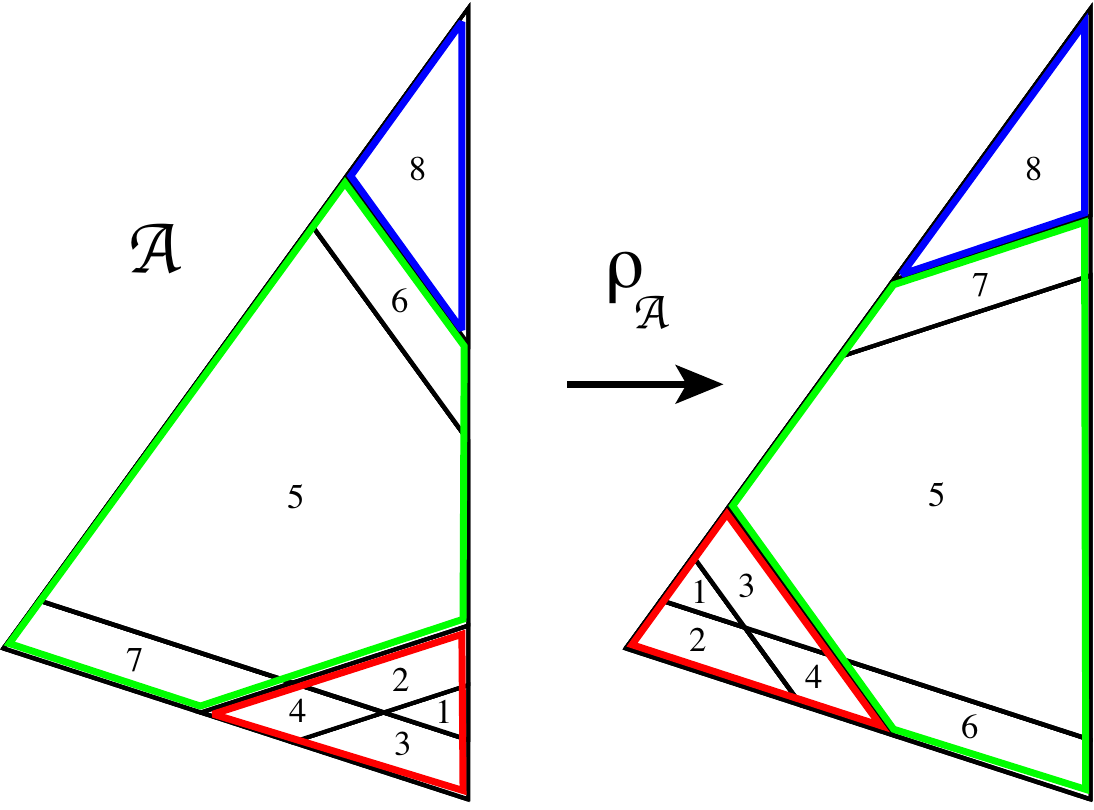,width=8 cm}\hfil

\caption{\label{fig:Apwi}\small 
The dressed domain $\mathcal{A}$, with its 8 atoms numbered as in (15).  
The boundaries of the composite atoms $C_1,C_2,C_3$ are coloured red, green, and blue, respectively.
}
\end{figure}

Our goal is to determine a parametric dressed domain $\mathcal{A}$ with 
bifurcation-free subdomain $\Pi(\mathcal{A})\subset\Pi(\mathcal{R})$ 
over which $\mathcal{R}$ is renormalizable.
To this end, we choose a specific parameter pair close to the $s=t$ boundary: $(s_0,t_0)=(-19/200,-1/10)$,
and we establish that at this value the renomalization is amenable to a
complete analysis (with computer assistance).
Specifically, we consider a three-step induction, first on the triangular atom $R_3$, 
followed by two inductions on sub-triangles, as shown in figure \ref{fig:3Xinduction}.
The last induction produces the dressed domain $\mathcal{A}$,
shown in figure \ref{fig:Apwi}, which is given by:
\begin{equation}\label{eq:A}
\mathcal{A}=(A,(A_1,\ldots,A_8), (\rho_{A_1},\ldots,\rho_{A_8}) ),
\end{equation}
with
\begin{eqnarray}
A&=&[ (2, 1, 0), (1, -1, 1), (-1 - s, \beta^4-s, 1 - s)]\nonumber\\
A_1&=&[(1, 0, 4), (-1, 1, -1), (\beta^4- t, 
  1 - s,\alpha \beta^4 - t)],\nonumber\\
A_2&=&[ (1, 4, 0, 4), (-1, 1, 1, -1), 
  (\beta^4 - t, \alpha \beta^4 - t,  1 - s, \alpha \beta^4- s)],\nonumber\\
A_3&=&[ (0, 1, 4, 1), (1, 1, -1, -1),
  (1 - s, \beta^4 - t, \alpha \beta^4 - t,  \beta^4 - s)],\nonumber\\
A_4&=&[ (4, 1, 4, 1), (1, 1, -1, -1),
  (\alpha \beta^4 - t, \beta^4 - t, \alpha \beta^4 - s,  \beta^4 - s)],\nonumber\\
A_5&=&[ (2, 1, 4, 0, 3), (1, -1, 1, 1, -1)\label{eq:Adef}\\
   && (-1 - s, \beta^4 - t, \alpha \beta^4 - s, 1 - s,  -4 \beta^3 - t)],\nonumber\\
A_6&=&[ (2, 3, 0, 3), (1, 1, 1, -1), (-1 - s, -4 \beta^3 - t, 1 - s, 
   -4 \beta^3 - s)],\nonumber\\
A_7&=&[ (2, 1, 4, 1), (1, -1, 1, 1),
   (-1 - s, \beta^4 - s, \alpha \beta^4 - s, \beta^4 - t)],\nonumber\\
 A_8&=&[ (2, 3, 0), (1, 1, 1), (-1 - s, -4 \beta^3 - s, 1 - s)],\nonumber
\end{eqnarray}
\begin{equation}\label{eq:rhoA}
\begin{array}{l}
\rho_{A_1}=\rho_{A_2}=\rho_{A_3}=\rho_{A_4}=\mathtt{T}_{(8 \beta^3, -2 \beta^5)}\, \mathtt{R}_2,\\
\rho_{A_5}=\rho_{A_6}=\rho_{A_7}=\mathtt{T}_{(2, 1+\beta^5)}\, \mathtt{R}_3,\\
\rho_{A_8}=\mathtt{T}_{(2-\beta^6, -\beta^5)}\, \mathtt{R}_2.
\end{array}
\end{equation}

\begin{figure}[h]
\hfil\epsfig{file=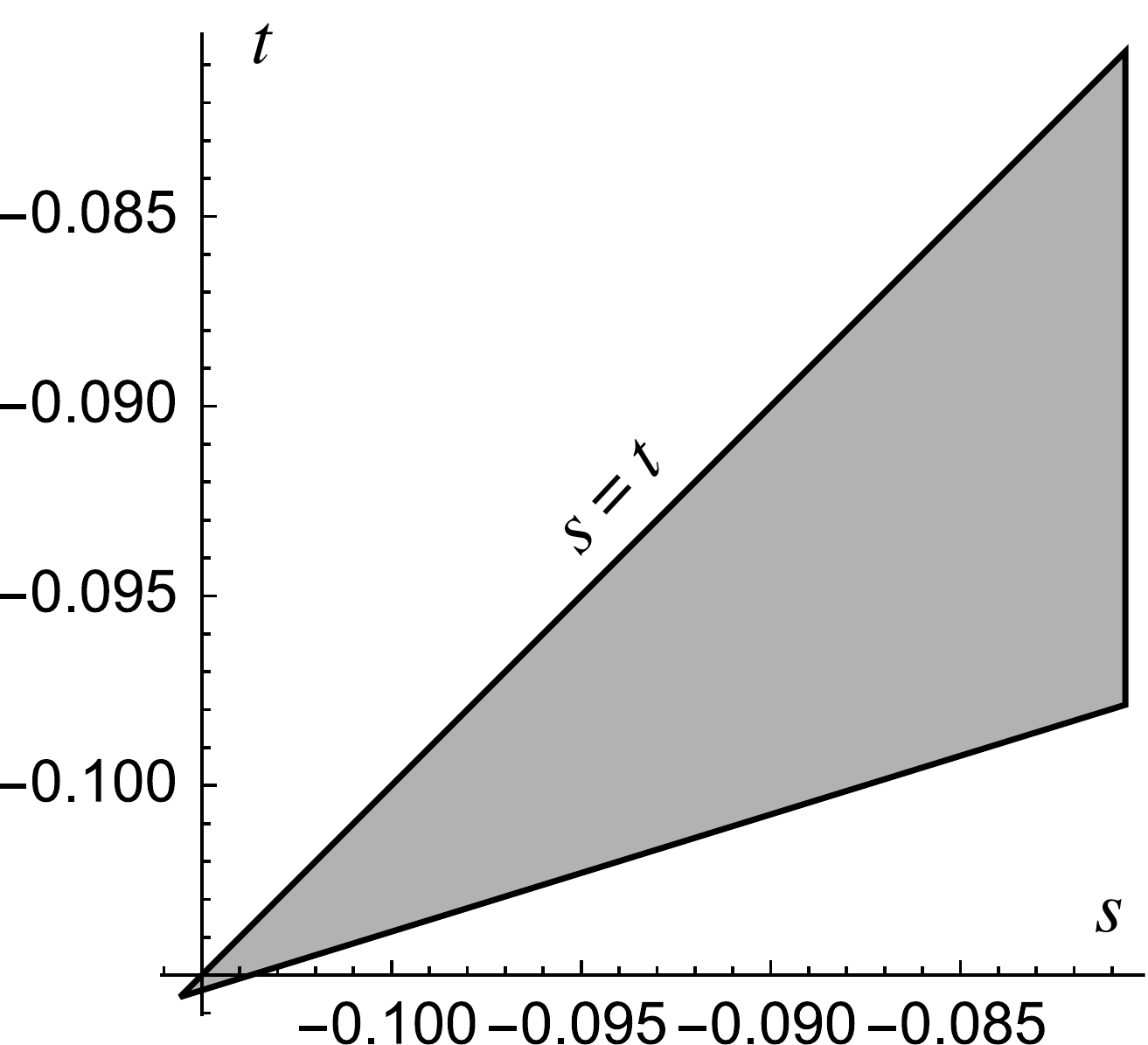,width=6 cm}\hfil
\caption{\label{fig:paramA}\small The parametric domain $\Pi(\mathcal{A})$.
}
\end{figure}

We find that $\Pi(\mathcal{A})$ is the triangle with vertices
$$
\displaystyle (-1+\frac{2}{\alpha},-1+\frac{2}{\alpha}),
\quad 
\displaystyle (\frac{1}{2}(11-5\alpha),\frac{1}{4}(-25+11\alpha)),
\quad
\displaystyle (\frac{1}{2}(11-5\alpha),\frac{1}{2}(11-5\alpha)),
$$
shown in figure \ref{fig:paramA}. 
One verifies that $\Pi(\mathcal{A})$ is adjacent to the line $s=t$ and that $(s_0,t_0)$
lies in its interior.

\begin{figure}[h]
\hfil\epsfig{file=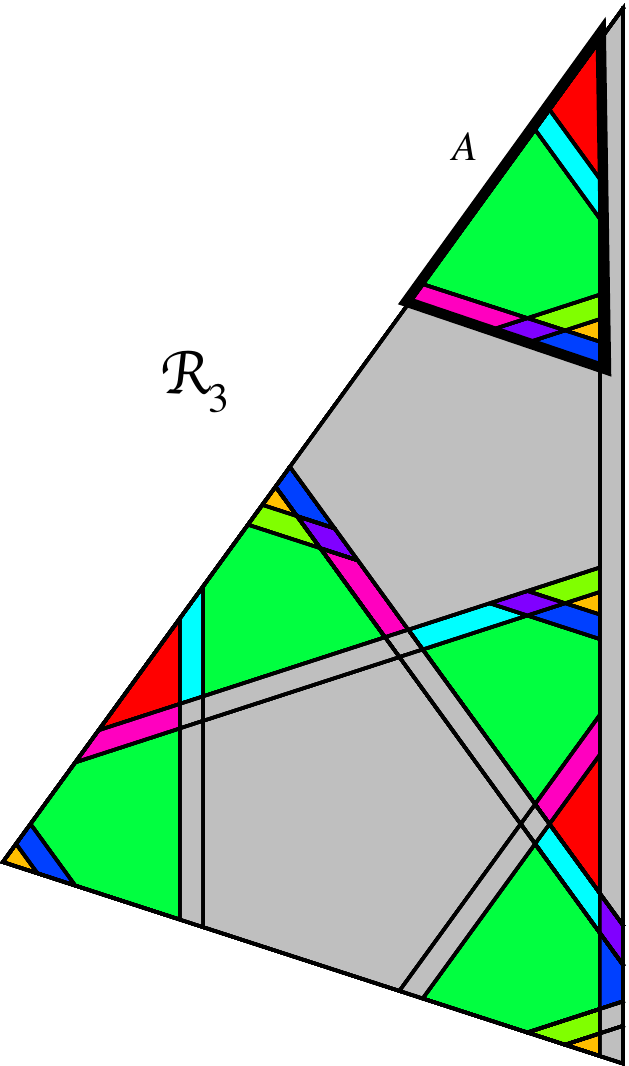,width=6cm}\hfil
\caption{\label{fig:tilingA}\small 
Tiling of $\mathcal{R}_3$ by return orbits of the 8 atoms of $\mathcal{A}$ (coloured) and 7 periodic tiles (grey).  
Note that  $A_1,A_2,A_3,A_4$, which comprise $C_1$, have distinct return orbits.
}
\end{figure}

Using direct iteration, we verify that for all $(s,t)\in\Pi(\mathcal{A})$ the return orbits of the eight 
atoms of $\mathcal{A}$, together with those of 13 periodic tiles, completely tile the rhombus $R$ 
(see figure \ref{fig:tilingA}).

A decisive simplification of the analysis results from the observation that
the four atoms $A_1,\ldots,A_4$ of $\mathcal{A}$ are mapped by the same isometry, and hence,
with regard to the first-return map to $\mathcal{A}$, can be merged into a single triangular tile, $A_{1234}$.
Similarly, atoms $A_5,A_6,A_7$ can be merged into a single reflection-symmetric pentagon, $A_{567}$. 
The mergers have been suggested in the shading of the tiles in figure \ref{fig:Apwi}.
The dressed domain thus simplifies into
\begin{eqnarray}
\mathcal{C}&=&\left(C, (C_1,C_2,C_3), (\rho_{C_1},\rho_{C_2},\rho_{C_3})\right)\nonumber\\
& \stackrel{\rm def}{=}& \left(A, (A_{8},A_{567},A_{1234}), (\rho_{A_8},\rho_{A_5},\rho_{A_1})\right) 
\label{eq:Cdef}
\end{eqnarray}
Moreover, one verifies that over $\Pi(\mathcal{C})=\Pi(\mathcal{A})$, we have
$\mathcal{C}\sim \mathcal{B}$.
The intrinsic shape parameter of $\mathcal{C}$ can be calculated from the ratio $\eta_\mathcal{C}$ 
of the altitude of $C_3$ to that of $C$:
\beq\label{eq:tauC}
\tau_\mathcal{C}=\frac{\omega^2 \eta_\mathcal{C}}{1+ \eta_\mathcal{C}} = \omega^7(\alpha s + \beta^3).
\eeq
As we transverse $\Pi(\mathcal{C})$ from left to right, $s$ increases from $2/\alpha-1$ to $(11-5\alpha)/2$, 
with $\tau_\mathcal{C}$ increasing from $0$ to $1$.

The issue of recursive tiling is now rather subtle. The rhombus is certainly tiled by the return
orbits of $C_1, C_2, C_3$, and the 13 periodic tiles which arose in the induction on $\mathcal{A}$. 
(see figure \ref{fig:tilingA}). 
However, the return paths are not the same for all tiles.
In (say) $C_3$, the tiles $A_1,A_2,A_3,A_4$ have four distinct 78-step return paths, which go their 
separate ways, but recombine eventually to form an atom of the dressed domain with a unique isometry.
The coincidence of the return times is not necessary to the recombination, as these times could differ by 
any integer multiples of 5. As a result, the partition of $C_3$ into $A_1,A_2,A_3,A_4$ is relevant 
to the recursive tiling of the original rhombus, but not to dynamical self-similarity.
(We shall encounter again the same phenomenon ---recombination with different
return times--- in section \ref{section:Explorations}.)

The parameter pairs $(s,t)\in\Pi(\mathcal{C})$ corresponding to self-similarity for the rhombus map $\mathcal{R}$ 
are now determined by the self-similarity of the induced dressed subdomain $\mathcal{C}$. 
In turn, the latter are the values of $s$ for which the base triangle is self-similar, 
namely $\tau_\mathcal{C}(s)\in (0,1)\cap \Q(\omega)$, while $t$ is unconstrained.
Thus, by (\ref{eq:tauC}), $\mathcal{R}$ is renormalizable in $\Pi(\mathcal{C})$ if and only if
$$
(s,t)\in \Pi(\mathcal{C}) \cap \left(\mathbb{Q}(\omega)\times \mathbb{R}\right).
$$ 
This is the main result of this section.

\clearpage
\section{Continuous self-similarity in Rauzy induction}\label{section:Rauzy}
Recombination of atoms, and the resulting appearance of a free parameter in self-similarity
may seem a coincidental feature of planar PWI's. 
This phenomenon is in fact common in the Rauzy-Veech analysis of renormalizable interval 
exchange transformations (IET's) \cite{Rauzy,Veech,Viana}.

We fix a half-open interval $\Omega=[0,l)$ and a partition of $\Omega$ into $n$ half-open 
sub intervals $\Omega_i$.
An IET is a piecewise isometry of $\Omega$ which is a translation on each $\Omega_i$.
We represent it as a pair $(\pi,\Lambda)$, where $\Lambda=(\lambda_1,\ldots,\lambda_n)$ is the vector of the lengths 
of the sub-intervals and $\pi$ is the permutation of $\{1,\ldots,n\}$ such that the
intervals in the image appear in the order $\pi(1),\ldots,\pi(n)$.

We assume that $\pi$ is {\it irreducible} in the sense that $\{1,\ldots,k\}$ is mapped into itself only if $k=n$.
If we fix $n$, then $(\pi,\Lambda)$ is a parametric PWI, with discrete and continuous parameters $\pi$ and 
$\Lambda$, respectively. 
IET's which differ only by an overall translation or scale transformation are considered equivalent.
 
The Rauzy-Veech induction on $(\pi,\Lambda)$ consists of inducing on the larger of the two intervals 
$\Omega^{(0)}=[0,l-\lambda_n)$ and $\Omega^{(1)}=[0,l-\lambda_{\pi^{-1}(n)})$, denoted by type 
0 and type 1 induction, respectively (the case $|\Omega^{(0)}|=|\Omega^{(1)}|$ is excluded from 
consideration, as in this case the map is not minimal). Induction corresponds to a map 
$(\pi,\Lambda)\mapsto (\pi',\Lambda')$.
Letting 
$$
(a_i(\pi), A_i(\pi)^{-1}\Lambda)=(\pi',\Lambda')\qquad i=0,1
$$
one finds that $A_i(\pi)^{-1}$ is an $n\times n$ integral matrix (see \cite{PoggiaspallaLowensteinVivaldi} 
for explicit expressions for $a_i$ and $A_i(\pi)$).

\begin{figure}[b]
\hfil\epsfig{file=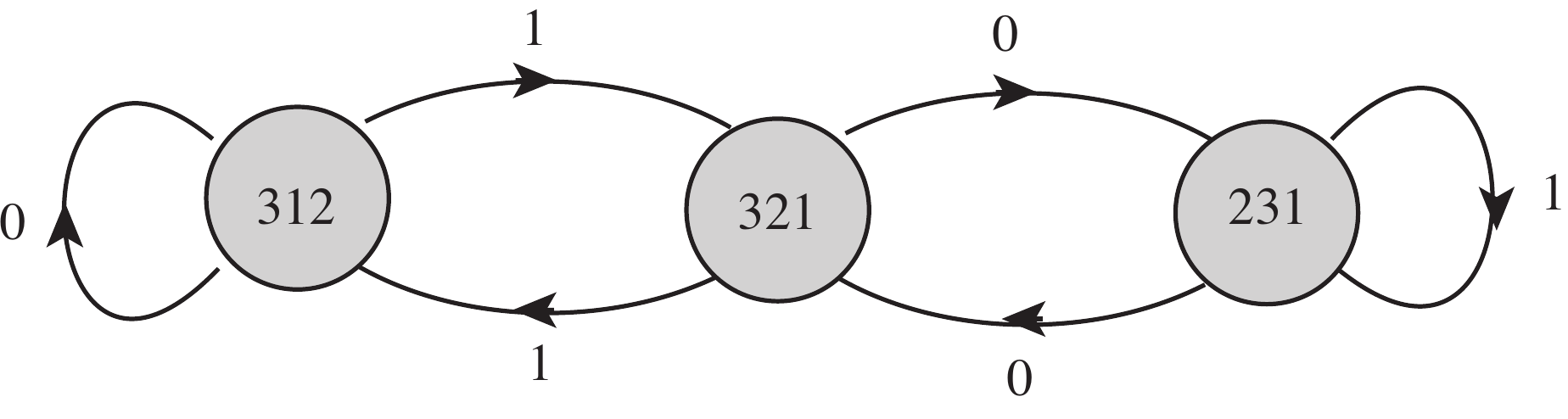,width=10cm}\hfil
\caption{\label{fig:Rauzygraph1}\small Rauzy graph for $n=3$.
All permutations are degenerate, and their translation surface is a torus.
Any renormalizable IET with three intervals will have a free parameter.
}
\end{figure}

\begin{figure}[h]
\hfil\epsfig{file=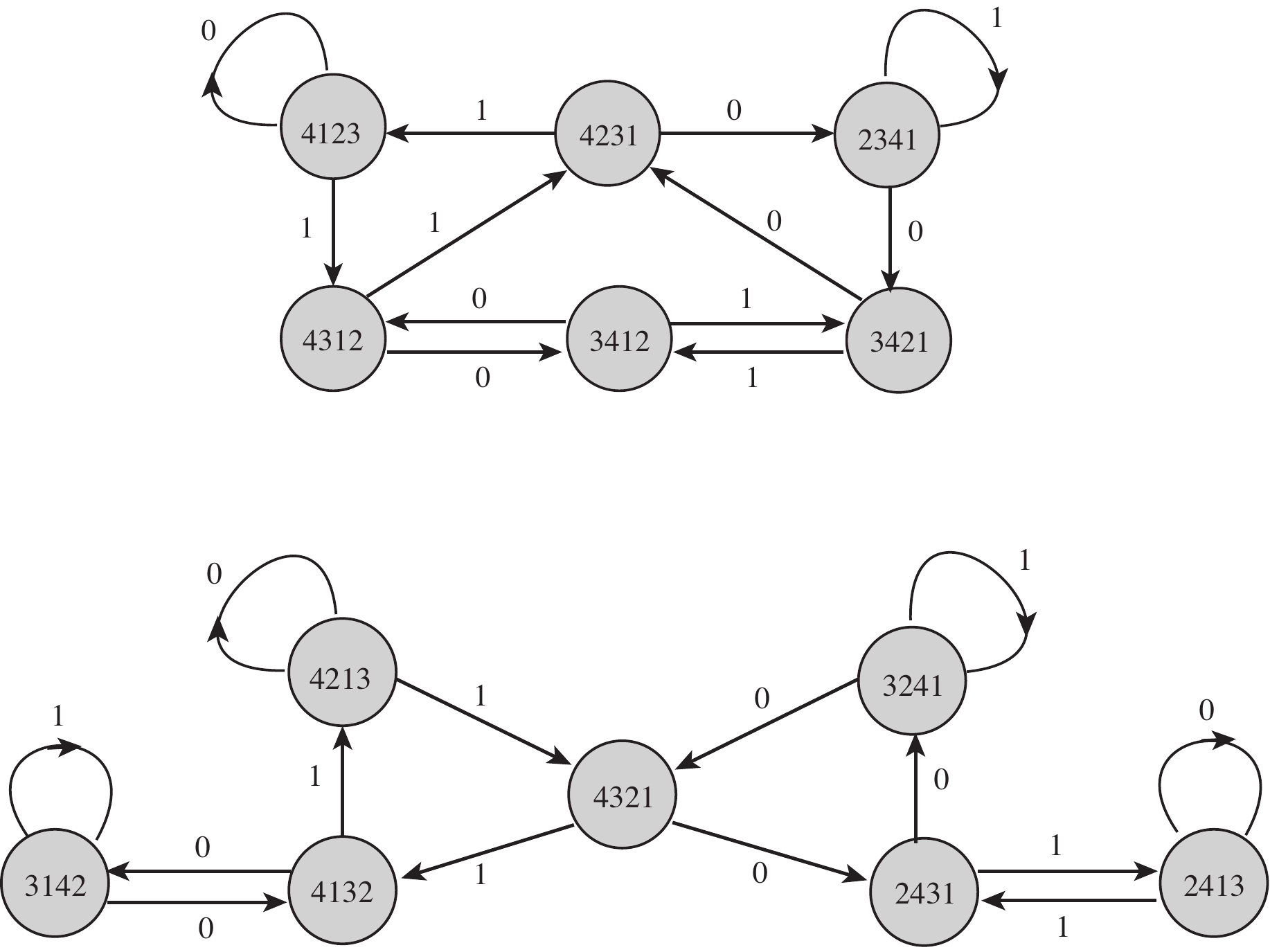,width=10cm}\hfil
\caption{\label{fig:Rauzygraph2}\small Rauzy graph for $n=4$. All permutations
of the upper component are degenerate. No permutation of the lower component has
this property.
}
\end{figure}

The permutations $\pi$ of $n$ symbols are then represented as the vertices of the 
{\it Rauzy graph}.
Each vertex has two outgoing and two incoming edges, associated with $a_i$ and 
$a_i^{-1}$, respectively, for $i=0,1$. 
The \textit{Rauzy classes} are the connected components of the graph.
These IET's are linked by a sequence of Rauzy inductions, and a self-similar 
IET corresponds to a path on a Rauzy graph which terminates in a closed circuit,
$\pi_1,\pi_2,\ldots,\pi_p$.
Transversing such a circuit produces an induced IET which is a rescaled version of 
the original one. Its length vector is an eigenvector of a product of matrices 
$A_{i_k}(\pi_k)^{-1}$, $k=1,\ldots,p$, with a scale factor given by the corresponding 
eigenvalue. In figures \ref{fig:Rauzygraph1} and \ref{fig:Rauzygraph2} we display 
the Rauzy graphs for $n=3$ and $n=4$ \cite[section 6]{Viana}. 

Two of the three irreducible permutations of the Rauzy graph for $n=3$ 
(figure \ref{fig:Rauzygraph1})
and all of the permutations for the first class for $n=4$ (figure \ref{fig:Rauzygraph2})
contain a consecutive pair: $(\ldots,j,j+1,\ldots)$. 
For those IET's, the consecutive intervals $\Omega_j$ and $\Omega_{j+1}$ have the same translation
vector, and hence the interval $\Omega_{j}\cup\Omega_{j+1}$ may be merged into a single interval
of length $\lambda'_j=\lambda_j+\lambda_{j+1}$. The $n$-interval IET is thus equivalent to an 
$n-1$-interval IET. If the latter is self-similar, then the original PWI is also self-similar 
for any choice of the parameter $\lambda_j\in[0,\lambda_j']$. 

Accordingly, we say that a permutation is \textit{degenerate} if it has consecutive pairs or if it acquires this 
property after a single induction. In the latter case, consecutive atoms of the child IET 
have distinct return paths in the tiling of the parent. 

Such a degeneracy is best understood by representing an IET as Poincar\'e section 
of a flow on a translation surface \cite{Veech,KontsevichZorich,Viana}.
The latter is a polygon with $2n$ sides ($n$ is the number of intervals), labelled according to the ordering of 
the intervals before and after the permutation. The sides that correspond to the same 
interval have equal length and are parallel, and they are to be identified. A rectilinear
flow on the plane will develop conical singularities on the surface, in correspondence 
to the vertices of the $2n$-gon.
While the translation surface is not unique, its genus and singularities (given 
by the total angle $2\pi(m+1)$ at the identified vertices) depend only on the Rauzy class. 
The removable singularities ($m=0$) correspond to degenerate permutations, and they signal 
the appearance of free parameters in renormalizability.
Since the translation surface does not change under induction, these structures 
depend only on the Rauzy class to which the permutation belongs.
Furthermore, for any $n\geqslant 4$
there are both degenerate and non-degenerate permutations,
e.g., $(n,\ldots,n-1,1)$ and $(n,n-1,\ldots,1)$, respectively.

The permutations which acquire consecutive pairs after induction are uniquely of the 
form $(\ldots, n, k+1,\ldots,k)$ for some $k<n-1$, and all four of its neighbours 
in a Rauzy graph have consecutive pairs, thanks to the relations
\begin{eqnarray*}
a_0\left((. . . n, k+1, . . ., k)\right) &=& (. . ., , k+1, k+2, . . .,k),\\
a_0^{-1}\left((. . . n, k+1, . . ., k)\right)&=& (. . ., n-1, n, . . .,k),\\
a_1\left((. . . n, k+1, . . ., k)\right) &=& (. . ., n, k, k+1, . . .),\\
a_1^{-1}\left((. . . n, k+1, . . ., k)\right) &=& (. . ., n, . . .,k, k+1).
\end{eqnarray*}
The $n=3$ class provides the simplest illustration of this phenomenon.

\clearpage
\section{Weakly-discrete self-similarity}\label{section:Explorations}
The analogy with Rauzy induction suggests that there might exist two-parameter planar 
PWI's which do not admit self-similarity with free parameters, meaning that both 
parameters would be algebraically constrained.
We shall exhibit a weak form of this property, resulting from the coexistence
of two systems with continuous self-similarity.

\begin{figure}[h]
\hfil\epsfig{file=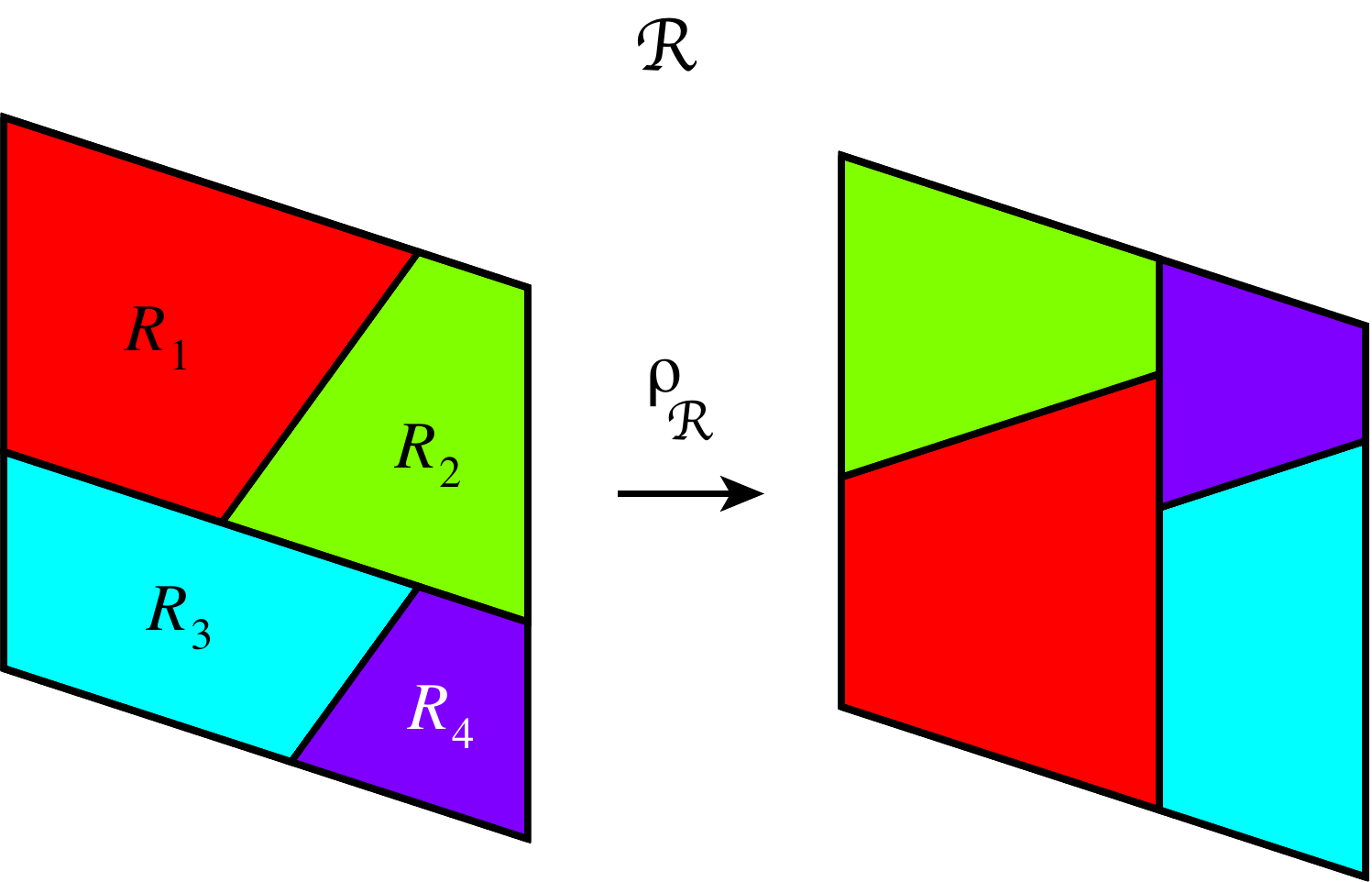,width=12 cm}\hfil
\caption{\label{fig:4Rpwi} \small Four-atom rhombus PWI.} 
\end{figure}
\begin{figure}[h]
\hfil\epsfig{file=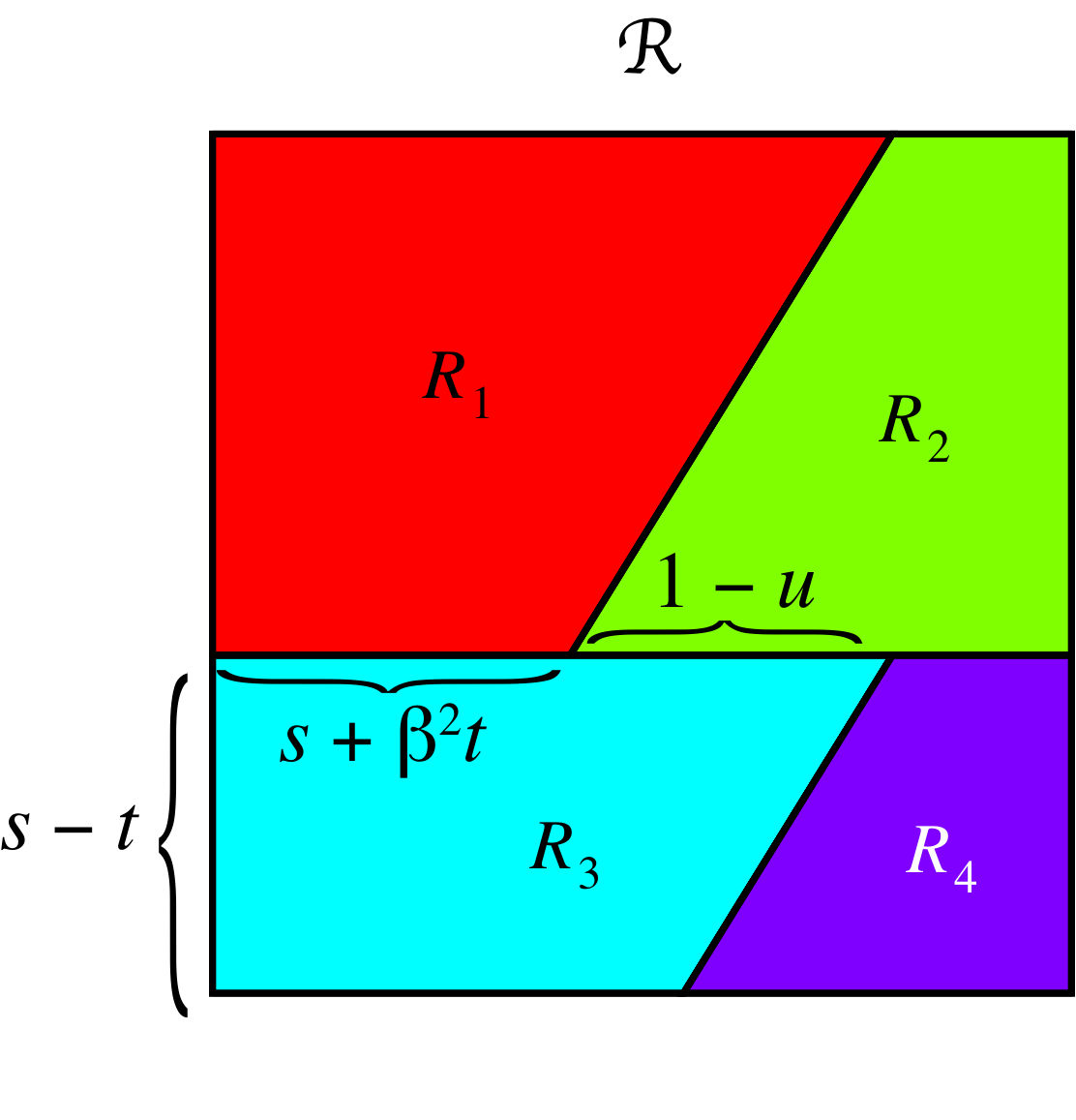,width=7 cm}\hfil
\caption{\label{fig:4square} \small Dressed domain $\mathcal{R}$. 
Dependence of atoms on parameters $s,t,u$ is shown.} 
\end{figure}
Our starting point is the four-atom, three-parameter PWI of the $2\pi/5$ rhombus 
shown in figure \ref{fig:4Rpwi}.  The interpretation of the parameters 
$s,t,u$ is made clear in the conjugate system of figure \ref{fig:4square}, where the 
rhombus appears as a unit square. 
(A similar strategy can be pursued for the five-atom family of figures
\ref{fig:5Rpwi} and \ref{fig:5square}, but we found that the four-atom family is
somewhat easier to work with.)

In convenient coordinates, the dressed domain is
$$
\mathcal{R}=(R,(R_1,\ldots,R_4), (\rho_1,\ldots,\rho_4) ),
$$
with
\begin{eqnarray}
R&=& [(0, 1, 0, 1), (-1, -1, 1, 1), (-t, -s, 1 - t, 1 - s)],\nonumber\\
R_1&=&[(0,1,2,1),(-1,-1,-1,1),(-t,-t,-s,1-s)],\nonumber\\
R_2&=&[(0,1,2,1),(1,1,1,-1),(1-t,1-s,-s,-t)],\nonumber\\
R_3&=&[(0,1,2,1),(-1,-1,-1,1),(-t,-s,-1-s+u,-t)],\label{eq:4Rdef}\\
R_4&=&[(0,1,2,1),(1,1,1,-1),(1-t,-t,-1-s+u,-s)]\nonumber
\end{eqnarray}
\begin{equation}\label{eq:rho4Rdef}
\begin{array}{lll}
\rho_1=\mathtt{T}_{(0,0)}\, \mathtt{R}_4,&\qquad
\rho_2=\mathtt{T}_{(0,1)}\, \mathtt{R}_4,\\
\rho_3=\mathtt{T}_{(1,1-u)}\, \mathtt{R}_4&\qquad
\rho_4=\mathtt{T}_{(1,2-u)}\, \mathtt{R}_4.
\end{array}
\end{equation}
From figure \ref{fig:4square}, we see that the bifurcation-free domain $\Pi(\mathcal{R})
\subset\mathbb{R}^3$ is the polytope bounded by the planes $s-t=0$, $s-t=1$, $u=1$, 
$s+\beta^2t=0$, $u-s-\beta^2 t=0$, $\beta^2-\beta^2 s-t=0$, and  $1+\beta^2 s+t-u=0$. 

This system has a simple one-parameter subsystem on the line $L$ defined by $s-t=\beta^2, u=\beta$.
We shall consider a two-parameter perturbation of this subsystem in the plane $u=\beta$ which intersects $L$.  
(We have also considered other planes, obtaining other manageable examples:
see remarks at the end of this section.)

\begin{figure}[h]
\hfil\epsfig{file=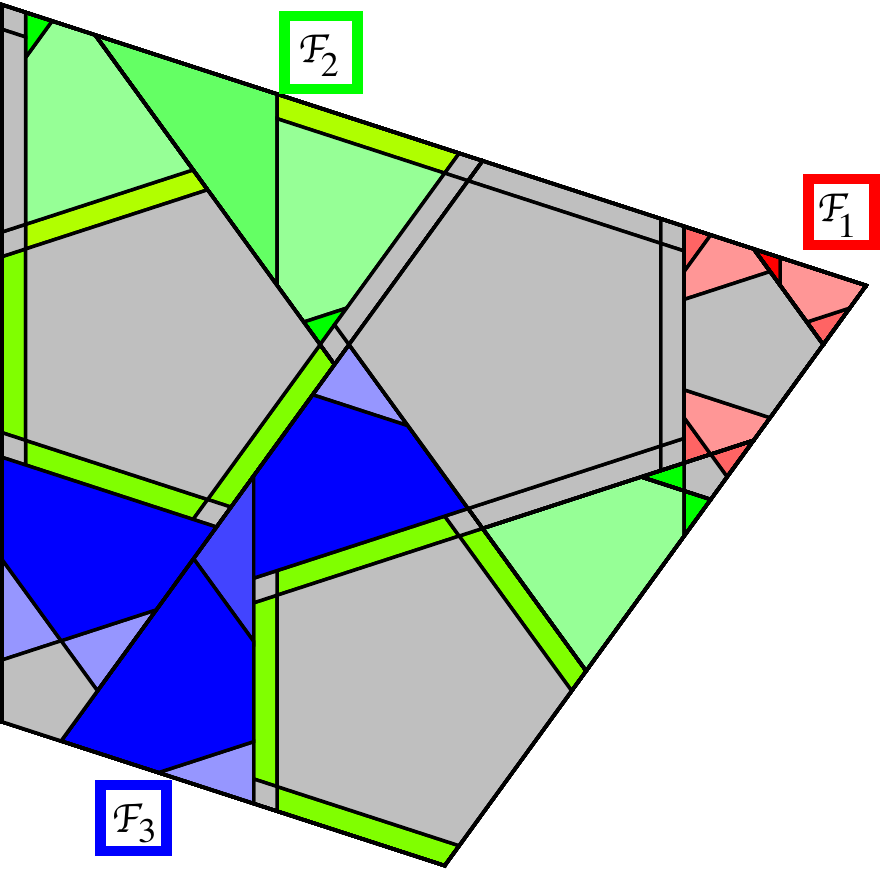,width=7 cm}\hfil
\caption{\label{fig:Ftiling} \small Tiling of $F$ by return orbits of 
$\mathcal{F}_1$ (red), $\mathcal{F}_2$ (green), $\mathcal{F}_3$ (blue), and seven periodic tiles (grey).} 
\end{figure}

\begin{figure}[h]
\hfil\epsfig{file=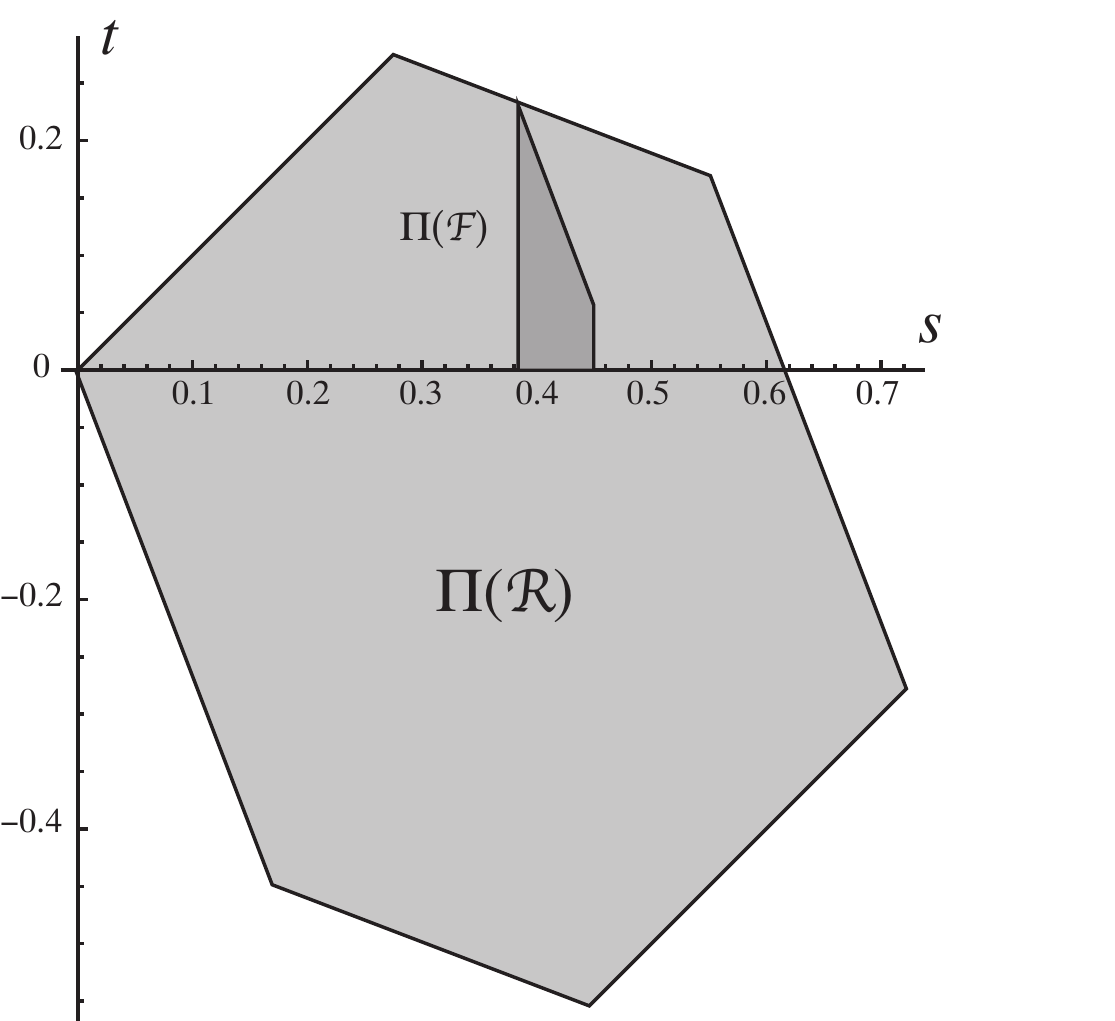,width=8 cm}\hfil
\quad
\hfil\epsfig{file=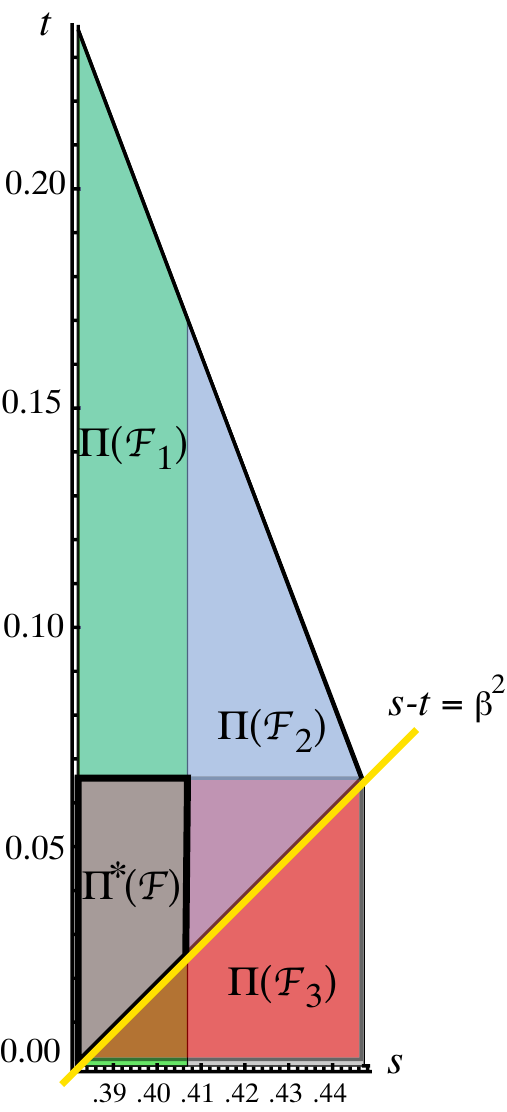,width=4 cm}\hfil
\caption{\label{fig:paramF} \small 
Left: The bifurcation-free domains $\Pi(\mathcal{R})$ and $\Pi(\mathcal{F})$ for the parameters $s$ and $t$, with $u=\beta$.
Right: Detailed view of $\Pi(\mathcal{F})$, showing the trapezoidal domain $\Pi^*(\mathcal{F})$ of
equation (\ref{eq:PiStarF}). The latter is constructed as the intersection of $\Pi(\mathcal{F}_1)$ 
(semi-transparent green trapezoid),  $\Pi(\mathcal{F}_2)$ (semi-transparent blue triangle), and 
$\Pi(\mathcal{F}_3)$ (semi-transparent red square).  
The seven periodic tiles $\mathcal{P}_i$ do not contribute any additional constraints. 
The unperturbed one-parameter model corresponds to the (yellow) line $s-t=\beta^2$, 
which lies along the south-east boundary of $\Pi(\mathcal{F})$ }.
\end{figure}

Setting $u=\beta$, the parameter polytope reduces to the hexagonal domain shown in figure \ref{fig:paramF} (left).
As done in section \ref{section:MainResults}, we choose a parameter pair close to $L$ lying within
such a domain: $(s_0,t_0)=(2/5,1/25)$.
\begin{figure}[h]
\hfil\epsfig{file=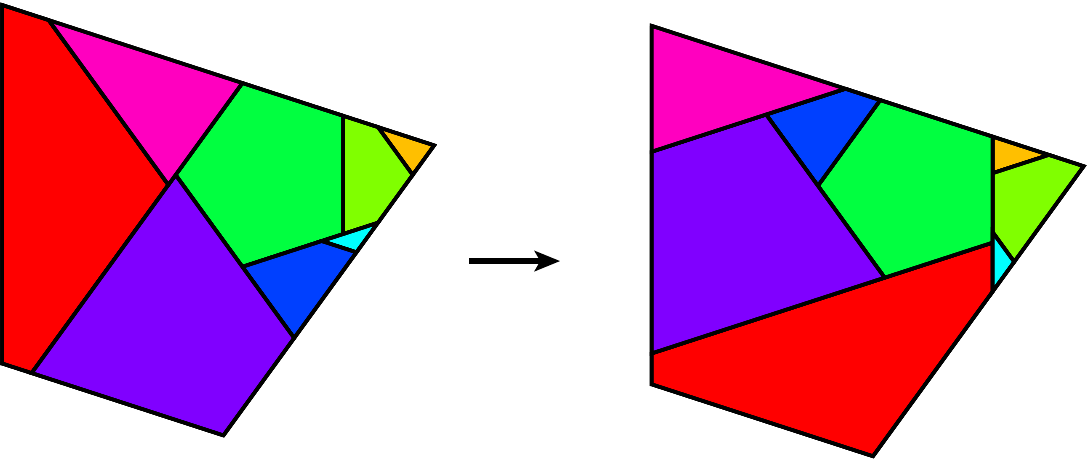,width=12 cm}\hfil
\caption{\label{fig:Fpwi} \small Induced dressed domain $\mathcal{F}$.} 
\end{figure}
By inducing on the trapezoidal atom $R_1$, we obtain the parametric dressed domain $\mathcal{F}$ shown 
in figure \ref{fig:Fpwi}. One readily verifies that the return orbits of the eight atoms of $\mathcal{F}$ 
completely tile $R$, so that the renormalizability of $\mathcal{R}$ will follow from that of $\mathcal{F}$. 
We find that the complete tiling of $F$ by renormalizable dressed sub-domains, given in figure \ref{fig:Ftiling}, requires the return 
orbits of three dressed triangles $\mathcal{F}_1, \mathcal{F}_2, \mathcal{F}_3$, plus 
seven periodic tiles $\mathcal{P}_i$ (five regular pentagons, one trapezoid, and one rhombus).  

Letting 
\begin{equation}\label{eq:PiStarF}
\Pi^*(\mathcal{F})=\bigcap_{i=1}^3 \Pi(\mathcal{F}_i)\,\bigcap_{i=1}^7\Pi(\mathcal{P}_i)
\end{equation}
we find that $\Pi^*(\mathcal{F})$ is the quadrilateral with vertices

$$
\left(\beta^2,0\right),\quad
\left(\beta^2+\beta^6/\alpha,\beta^6/\alpha)\right),\quad
\left(\beta^2+\beta^6/\alpha,\beta^4/\alpha\right),\quad
\left(\beta^2,\beta^4/\alpha\right),
$$
shown in figure \ref{fig:paramF} (right).   
Note that one of the bounding edges of the parameter domain coincides with the 
line $L$ which was the starting point of our perturbative exploration.
 
The dressed domains $\mathcal{F}_1$ and $\mathcal{F}_2$ are base triangles equivalent to the 
prototype $\mathcal{B}$, with respective shape parameters 
$$
\tau_1=\omega^6 \alpha (s-\beta^2),\qquad
\tau_2=\omega^4 \alpha t.
$$
Examination of $\mathcal{F}_3$ shows that its atoms with four and five sides share the same 
isometry (in spite of having different return paths, and even different return times on the rhombus), 
and hence can be merged for the purpose of testing renormalizability.  
After the merger, $\mathcal{F}_3$ is also equivalent to $\mathcal{B}$, with shape parameter 
$$
\tau_3=\omega^6 \alpha (s-\beta^2).
$$

Since each  $\mathcal{F}_i$ is renormalizable for $\tau_i\in \mathbb{Q}(\sqrt{5})$, 
we conclude that $\mathcal{F}$ (hence $\mathcal{R}$) is renormalizable when all three 
shape parameters are in $\mathbb{Q}(\sqrt{5})$, i.e., when $(s,t)$ is 
constrained to belong to $\Pi^*(\mathcal{F})\cap\mathbb{Q}(\sqrt{5})^2$. 

Each of the three renormalizable dressed domains $\mathcal{F}_i$ provides a sequence of nested coverings 
of a distinct invariant component of the exceptional set complementary 
to all periodic points of the rhombus.  
The number of distinct ergodic components of the exceptional set is thus at least three.   
For the model of section \ref{section:MainResults}, on the other hand, we believe
that there is a single ergodic component.

In closing, we summarize briefly the results of our explorations of the
three-parameter space of the four- or five-atom rhombus maps.

Manageable renormalizations are likely to be found in systems specified by parameters
of small height. [The height $H(\zeta)$ of the algebraic number $\zeta=(m/n)+(m'/n')\omega$ 
is defined as $H(\zeta)=\max(|m|,|n|,|m'|,|n'|)$].
Such was the case for the domain $\mathcal{R}$, and the two-parameter restriction $u=\beta$
described above. We have considered other planes of small height: $s-t-\omega u+\beta=0$, 
$s-t+u-1=0$, etc.
Within such planes, and for carefully chosen parameter patches, we have encountered a number 
of two-parameter renormalizable models. 
These are characterized by a decomposition of the rhombus into $N$ disjoint dressed domains,
each tiled by the return orbits of a single base triangle (provided we ignore the 
``decorations'' produced by the common boundaries of merged atoms), plus a finite number of 
periodic tiles.

For the plane $u=0$ of the five-atom map, the case $N=1$ appears to be the norm,
so that in those models the exceptional set is likely to be uniquely ergodic.
Elsewhere, a proliferation of ergodic components is typical (albeit not universal).
Notably, we have found no example of a rigidly self-similar, single component 
piecewise isometry with two or more parameters.
Whether such a dynamical system exists at all remains an important open question.
%

\end{document}